\documentclass[12pt]{article}

\usepackage[utf8]{inputenc}
\usepackage[a4paper,nomarginpar]{geometry}
\usepackage[english]{babel}
\usepackage{enumitem}
\usepackage{amstext,amsmath,amsfonts,amssymb,amsthm,amscd}
\usepackage{tikz,xcolor}
\usepackage{soul}
\usepackage{comment}
\usepackage{textcomp}
\usepackage{enumitem}
\usepackage{theoremref}
\usepackage{mathtools}
\usepackage{float}
\usepackage{etoc}
\usepackage{stmaryrd}

\allowdisplaybreaks

\theoremstyle{plain}
\newtheorem{theorem}{Theorem}
\newtheorem{corollary}[theorem]{Corollary}

\theoremstyle{definition}
\newtheorem{definition}[theorem]{Definition}
\newtheorem{remark}[theorem]{Remark}
\newtheorem{example}[theorem]{Example}

\newcommand{\N}{\mathbb{N}} 
\newcommand{\Z}{\mathbb{Z}} 
\newcommand{\R}{\mathbb{R}} 
\newcommand{\C}{\mathbb{C}} 
\newcommand{\D}{\mathbb{D}} 
\newcommand{\T}{\mathbb{T}} 
\newcommand{\K}{\mathbb{K}} 

\newcommand{\cB}{\mathcal{B}}

\newcommand{\cP}{\mathcal{P}}

\newcommand{\cX}{\mathcal{X}}

\newcommand{\fM}{\mathfrak{M}}

\newcommand{\re}{\operatorname{Re}}

\newcommand{\card}{\operatorname{card}}

\newcommand{\udens}{\operatorname{\overline{dens}}}

\newcommand{\spa}{\operatorname{span}}      
\newcommand{\ran}{\operatorname{ran}}

\newcommand{\eps}{\varepsilon}
\newcommand{\ov}{\overline}
\newcommand{\wt}{\widetilde}

\title{Shadowing property  for  linear operators on Banach spaces}

\author{Nilson C. Bernardes Jr., Antonio Bonilla}

\date{}

\begin{document}

\maketitle

\begin{abstract}
This note is dedicated to the study of the shadowing property for (bounded linear) operators on Banach spaces. 
We combine the spectral characterization of the shadowing property, recently established by Dragi\v{c}evi\'c and Pituk, 
with tools from local spectral theory, to obtain new results and simpler proofs of known results on shadowing for operators. 
We also prove that the shadowing property always implies the average shadowing property 
and the asymptotic average shadowing property for invertible operators on Banach spaces.
Furthermore, we investigate the shadowing property for some specific classes of operators, 
such as weighted composition operators and adjoints of multiplication operators on certain spaces of holomorphic functions.
\end{abstract}

\let\thefootnote\relax

\begin{center}
\etocsetnexttocdepth{1} 
\etocsettocstyle{\small}{}
\tableofcontents
\end{center}


\section{Introduction}

The shadowing property is one of the fundamental concepts in the theory of dynamical systems and its applications. 
It originated from the work of Anosov, Bowen and Sina$\breve{\text{\i}}$ in the late 1960s and early 1970s 
and has been extensively investigated by several authors in different contexts. 
General references on this concept include the books \cite{AKatBHas95,KPal00,SPil99}.
In the setting of linear dynamics, it has long been known that invertible hyperbolic operators on Banach spaces
have the shadowing property \cite{AMor81,JOmb94} and that converse holds for invertible operators on finite-dimensional spaces
\cite{AMor81,JOmb94} and for invertible normal operators on Hilbert spaces \cite{MMaz00}.
The first examples of non-hyperbolic operators with the shadowing property were obtained by Bernardes et al.\ \cite{BCDMP},
thereby solving a long-standing open problem. 
Since the publication of the article \cite{BCDMP} in 2018, a systematic study of the shadowing property 
in the context of linear dynamics has been conducted by various authors. 
Among the many interesting results obtained, let us mention the following: 
\begin{itemize}
\item The equality “hyperbolicity = expansivity + shadowing” for invertible operators on Banach spaces \cite{NBerAMes21,CirGolPuj21}.
\item A complete characterization of the shadowing property for weighted shifts on classical Banach sequence spaces \cite{NBerAMes21}.
\item The emergence of the concept of generalized hyperbolicity \cite{BCDMP,CirGolPuj21}.
\item The structural stability of every generalized hyperbolic operator and the generalized Grobman-Hartman theorem \cite{NBerAMes20}
  (see also \cite{B}).
\item The equivalence between shadowing and finite shadowing for operators on Banach spaces and
  the failure of this equivalence for operators on Fr\'echet spaces \cite{BP}.
\end{itemize}

Quite recently, an important breakthrough was achieved by Dragi\v{c}evi\'c and Pituk \cite{DP}, 
who obtained a spectral characterization of the shadowing property for operators on Banach spaces.

In this note, we combine the spectral characterization of the shadowing property with tools from local spectral theory to establish 
new results---as well as simpler proofs of known results---concerning the shadowing property for operators on Banach spaces. 
Furthermore, we prove that the shadowing property implies the average shadowing property
and the asymptotic average shadowing property, and we investigate the shadowing property for specific classes of operators, 
including weighted composition operators and adjoints of multiplication operators on certain spaces of holomorphic functions.

The note is organized as follows. 

In Section~\ref{Section2}, we present some preliminaries concerning various types of spectra and the shadowing property for operators.

In Section~\ref{Section3}, we establish a characterization of the shadowing property involving local spectral theory, 
present several classes of operators that have the shadowing property, 
and show various contexts in which the concepts of shadowing and hyperbolicity coincide. 
Furthermore, we study several properties of operators that have the shadowing property, including properties related to chaotic behaviors.

In Section~\ref{Section4}, we present a list of classes of operators and some specific operators that do not have the shadowing property.

In Section~\ref{Section5}, we prove that every operator with the shadowing property also has the average shadowing property and the asymptotic average shadowing property.

In Section~\ref{Section6}, we present characterizations of the shadowing property for invertible weighted composition operators 
on certain Banach spaces of holomorphic functions on the unit disk, including Hardy spaces and weighted Bergman spaces.

In Section~\ref{Section7}, we present simple proofs for the known characterizations of the shadowing property for unilateral and bilateral weighted backward shifts on the classical Banach sequence spaces $\ell^p$ ($1 < p <\infty$) and $c_0$, as well as some additional results.

In Section~\ref{Section8}, we study the shadowing property for operators of the form $\phi(T)$ and their adjoints, where $T$ is 
a non-invertible isometry and $\phi$ is a non-constant holomorphic function defined on a neighborhood of the closed unit disk. 
The general result is applied to the special case of multiplication Toeplitz operators $M_{\phi}$ on the Hardy space $H^2(\D)$.

In Section~\ref{Section9}, we study the shadowing property for weighted composition operators $\psi C_{\varphi}$ 
on the Hilbert space $L^2[0,1]$ that are determined by an $L^2$-admissible pair $(\varphi,\psi)$ of mappings.

Finally, in Section~\ref{Section10}, we present sufficient conditions for the shadowing property 
for weighted composition operators on $L^p(\mu)$ spaces.


\section{Preliminaries on shadowing for operators}\label{Section2}

Throughout this note, $X$ denotes an arbitrary complex Banach space, unless otherwise indicated.
By an {\em operator} on $X$, we mean a bounded linear map $T : X \to X$. 
We denote by $L(X)$ the Banach algebra of all operators on $X$ and by $GL(X)$ the group of all invertible operators on $X$.
Recall that the {\em spectrum}, the {\em point spectrum} and the {\em surjective spectrum} of $T$ are defined by
\begin{align*}
\sigma(T) &= \{\lambda \in \C : \lambda I - T \mbox{ is not invertible}\},\\
\sigma_p(T) &= \{\lambda \in \C : \lambda I - T \mbox{ is not injective}\},\\
\sigma_s(T) &= \{\lambda \in \C : \lambda I - T \mbox{ is not surjective}\},
\end{align*}
respectively. 
Recall also that the {\em approximate point spectrum} $\sigma_a(T)$ of $T$ is the set of all $\lambda \in \C$ 
for which there is a sequence $(x_n)_{n \in \N}$ of unit vectors in $X$ such that $\|\lambda x_n - Tx_n\| \to 0$ as $n \to \infty$.
Moreover, the {\em spectral radius} of $T$ is defined by $r(T) = \max_{\lambda \in \sigma(T)} |\lambda|$.
The {\em spectral radius formula} asserts that
\[
r(T) = \lim_{n \to \infty} \|T^n\|^{\frac{1}{n}} = \inf_{n \in \N} \|T^n\|^{\frac{1}{n}}.
\]

An operator $T \in L(X)$ is said to have the {\em positive shadowing property} if, for every $\eps > 0$,
there exists $\delta > 0$ such that for any sequence $(x_n)_{n \in \N_0}$ satisfying
\[
\|x_{n+1} - Tx_n\| \leq \delta \ \ \text{ for all } n \in \N_0,
\]
called a {\em $\delta$-pseudotrajectory} of $T$, there exists a point $x \in X$ whose true trajectory 
{\em $\eps$-shadows} $(x_n)_{n \in \N_0}$, that is,
\[
\|x_n - T^nx\| < \eps \ \ \text{ for all } n \in \N_0.
\]
For $T \in GL(X)$, the {\em shadowing property} is defined by replacing $\N_0$ with $\Z$ in the definition above.

An operator $T \in GL(X)$ is said to be {\em hyperbolic} if $\sigma(T) \cap \T = \emptyset$, 
where $\T$ denotes the unit circle in the complex plane.
It follows from the classical Riesz decomposition theorem that $T$ is hyperbolic if and only if $X$ admits a direct sum decomposition
$X = M \oplus N$, where $M$ and $N$ are closed $T$-invariant subspaces of $X$ such that $r(T|_M) < 1$ and $r(T^{-1}|_N) < 1$.
It has long been known that every hyperbolic operator $T \in GL(X)$ has the shadowing property.
The following generalization of this result was obtained by Bernardes et al.\ \cite{BCDMP}.

\begin{theorem}[{\cite[Theorem~A]{BCDMP}}]
{ \it Let $T \in GL(X)$. If X admits a direct sum decomposition $X = M \oplus N$, where $M$ and $N$ are closed subspaces of $X$ 
such that $T(M) \subset  M$, $T^{-1}(N) \subset N$, $r(T|_M) < 1$ and $r(T^{-1}|_N) < 1$, then $T$ has the shadowing property.}
\end{theorem}

This result made it possible to demonstrate the existence of non-hyperbolic operators with the shadowing property
\cite[Theorem~B]{BCDMP}. 
An operator satisfying the property of the result above was termed {\em generalized hyperbolic} in \cite{CirGolPuj21}, 
where additional properties of these operators were investigated.

The following spectral characterizations of the shadowing property for operators were recently obtained 
by Dragi\v{c}evi\'c and Pituk \cite{DP}.

\begin{theorem}[{\cite[Theorem~2.2]{DP}}]\label{Spectral1}
An operator $T \in GL(X)$ has the shadowing property if and only if $\sigma_s(T) \cap \T = \emptyset$.
\end{theorem}

\begin{theorem}[{\cite[Theorem~2.4]{DP}}]\label{Spectral2}
An operator $T \in L(X)$ has the positive shadowing property if and only if $\sigma_s(T) \cap \T = \emptyset$.
\end{theorem}

Note that, by using the fact that the concepts of shadowing and positive shadowing coincide for invertible operators on Banach spaces 
\cite[Theorem~1]{BP}, we see that Theorem~\ref{Spectral1} follows immediately from Theorem~\ref{Spectral2}.


\section {Operators with the shadowing property}\label{Section3}

Given an operator $T \in L(X)$ and a closed subset $F$ of $\C$, 
the {\em glocal spectral subspace} $\cX_T(F)$ consists of all $x \in X$ for which there exists 
an analytic function $f : \C \backslash F \to X$ that satisfies 
\[
(T-\lambda I)f(\lambda) = x \ \text{ for all } \lambda \in \C \backslash F.
\]
For each $x\in X$, the {\em local spectral radius} of $T$ at $x$ is defined by 
\[
r_T(x) = \limsup_{n \to \infty} \|T^nx\|^{\frac{1}{n}}.
\]
For each $r > 0$, let $\D_r$ (resp.\ $\ov{\D}_r$) denote the open (resp.\ closed) disk with center $0$ and radius $r$ in the complex plane.
It is common to denote $\D_1$ simply by $\D$.

\begin{theorem}\label{glocal} 
For every $T \in GL(X)$, the following assertions are equivalent:
\begin{itemize}
\item [\rm (i)] $T$ has the shadowing property.
\item [\rm (ii)] There exists $r \in (0,1)$ such that
  \[
  X = \cX_T(\ov{\D}_r) + \cX_{T^{-1}}(\ov{\D}_r).
  \]
\item [\rm (iii)] There exists $r \in (0,1)$ such that 
  \[
  X= \{x \in X : r_T(x) \leq r\} + \{x \in X : r_{T^{-1}}(x) \leq r\}.
  \]
\end{itemize}
\end{theorem}

\begin{proof}
By Theorem \ref{Spectral1}, $T$ has the shadowing property if and only if
\[
\sigma_s(T) \cap \T = \emptyset,
\]
which is equivalent to the existence of a number $r \in (0,1)$ such that
\[
\sigma_s(T) \subset \ov{\D}_r \cup (\C \backslash \D_{\frac{1}{r}}).
\]
By \cite[Theorem 3.3.12(a)]{LN} and \cite[Proposition 3.3.1(f)]{LN}, the above inclusion holds if and only if
\[
X = \cX_T(\ov{\D}_r \cup (\C \backslash \D_{\frac{1}{r}})) = \cX_T(\ov{\D}_r) + \cX_T(\C \backslash \D_{\frac{1}{r}}).
\]
Since $T$ is invertible, \cite[Theorem 3.3.6]{LN} assures that the above equality can be rewritten as 
\[
X = \cX_T(\ov{\D}_r) + \cX_{T^{-1}}(\ov{\D}_r).
\]
This shows that (i) $\Leftrightarrow$ (ii).

The equivalence (ii) $\Leftrightarrow$ (iii) follows immediately from \cite[Proposition 3.3.13]{LN}.
\end{proof}

As a consequence of Theorem \ref{glocal}, we can obtain simple proofs of \cite[Theorem~A]{BCDMP} and \cite[Proposition~2.5]{P}.

\begin{corollary}[{\cite[Theorem A]{BCDMP}}]\label{GH-SP}
Every generalized hyperbolic operator $T \in GL(X)$ has the shadowing property.
\end{corollary}

\begin{proof}
Since $T$ is generalized hyperbolic, there exist closed subspaces $M$ and $N$ of $X$ such that
$X = M \oplus N$, $T(M) \subset M$, $T^{-1}(N) \subset N$, $r(T|_M) < 1$ and $r(T^{-1}|_N) < 1$.
Choose a real number $r$ such that $\max\{r(T|_M),r(T^{-1}|_N) < r < 1$. Then,
\[
M \subset \{x \in X : r_T(x) \leq r\} \ \ \text{ and } \ \ N \subset \{x \in X : r_{T^{-1}}(x) \leq r\}.
\]
Thus, condition (iii) in Theorem~\ref{glocal} holds and $T$ has the shadowing property.
\end{proof}

\begin{corollary}[{\cite[Proposition 2.5]{P}}] \label{Cor-Pseudo}
Let $T\in  GL(X)$. If there exists an operator $B : X \to X$ such that
\[
\limsup_{n \to \infty} \|T^nB\|^{\frac{1}{n}} < 1 \ \ \text{ and } \ \ \limsup _{n \to \infty} \|T^{-n}(I-B)\|^{\frac{1}{n}} < 1,
\]
then $T$ has the shadowing property.
\end{corollary}

\begin{proof}
The hypothesis implies the existence of a number $r \in (0,1)$ such that  
\[
B(X) \subset \{x \in X : r_T(x) \leq r\} \ \ \text{ and } \ \ (I-B)(X) \subset \{x \in X :  r_{T^{-1}}(x) \leq r\}.
\]
Since $X = B(X) + (I-B)X$, Theorem~\ref{glocal} shows that $T$ has the shadowing property.
\end{proof}

Given any mapping $B : X \to X$ (not necessarily continuous or linear), let us adopt the following notation:
\[
\|B\| = \sup\{\|B(x)\| : \|x\| \leq 1\},
\]
where the supremum can be equal to $\infty$.
Upon examining the previous proof more carefully, we observe that the only property of $B$ that we need is that 
$B(tx) = t B(x)$ for all $x \in X$ and all $t \geq 0$.
Therefore, we have the following generalization of the previous corollary.

\begin{corollary}\label{Cor-G-Pseudo}
Let $T\in  GL(X)$. If there exists a mapping $B : X \to X$ (not necessarily continuous or linear) such that
$B(t x) = t B(x)$ for all $x \in X$ and all $t \geq 0$,
\[
\limsup_{n \to \infty} \|T^n \circ B\|^{\frac{1}{n}} < 1 \ \ \text{ and } \ \ \limsup _{n \to \infty} \|T^{-n} \circ (I-B)\|^{\frac{1}{n}} < 1,
\]
then $T$ has the shadowing property.
\end{corollary}

An operator $T \in GL(X)$ that satisfies the hypothesis of Corollary~\ref{Cor-Pseudo} was termed {\em pseudo-hyperbolic} in \cite{MNST}. This class was further generalized in the recent preprint \cite{DP2}:
$T \in GL(X)$ is said to be {\em generalized pseudo-hyperbolic} if there exists a continuous mapping $B : X \to X$ 
which is complex-homogeneous (i.e., $B(\lambda x) = \lambda B(x)$ for all $x \in X$ and all $\lambda \in \C$) and satisfies
\[
\limsup_{n \to \infty} \|T^n \circ B\|^{\frac{1}{n}} < 1 \ \ \text{ and } \ \ \limsup _{n \to \infty} \|T^{-n} \circ (I-B)\|^{\frac{1}{n}} < 1.
\]
Corollary~\ref{Cor-G-Pseudo} implies that every generalized pseudo-hyperbolic operator $T \in GL(X)$ has the shadowing property.

\begin{remark}
Condition (ii) (or (iii)) in Theorem~\ref{glocal} shows that the shadowing property appears to be, 
in a sense, close to generalized hyperbolicity.
More precisely, by considering $M = \cX_T(\ov{\D}_r)$ and $N = \cX_{T^{-1}}(\ov{\D}_r)$, 
we have that $M$ and $N$ are subspaces of $X$ such that $X = M + N$, $T(M) \subset M$, $T^{-1}(N) \subset N$,
$r_T(x) \leq r < 1$ for all $x \in M$, and $r_{T^{-1}}(x) \leq r < 1$ for all $x \in N$.
Furthermore, if $M$ (resp.\ $N$) is closed, then $r(T|_M) \leq r < 1$ (resp.\ $r(T^{-1}|_N) \leq r < 1$) by \cite[Proposition~3.3.15]{LN}.
However, it was recently announced in the preprints \cite{MNST,P2} that there exist operators on certain Banach spaces 
that possess the shadowing property and are not generalized hyperbolic; 
but it was also announced in \cite{P2} that shadowing and generalized hyperbolicity 
actually coincide in the case of complex separable Hilbert spaces. 
Furthermore, it was recently announced in \cite{DP2} that the shadowing property is, in fact, 
equivalent to generalized pseudo-hyperbolicity for invertible operators on Banach spaces.
\end{remark}

The following known facts can also be trivially deduced from Theorem~\ref{glocal}.

\begin{corollary} For any $T\in GL(X)$ and any $\lambda \in \C$ with $|\lambda| = 1$, we have that:
\begin{itemize}
\item [\rm (a)] $T^{-1}$ has the shadowing property if and only if so does $T$. 
\item [\rm (b)] $\lambda T$ has the shadowing property if and only if so does $T$.
\end{itemize}
\end{corollary}

Let us now present a version of Theorem~\ref{glocal} for operators that are not necessarily invertible.

\begin{theorem}\label{glocal2} 
For every $T \in L(X)$, the following assertions are equivalent:
\begin{itemize}
\item [\rm (i)] $T$ has  the positive  shadowing property.
\item [\rm (ii)] There exists $r \in (0,1)$ such that $ X = \cX_T(\ov{\D}_r) + \cX_T(\C \backslash \D_{\frac{1}{r}})$.
\item [\rm (iii)] There exists $r \in (0,1)$ such that $ X = \{x \in X : r_T(x) \leq r\} + \cX_T(\C \backslash \D_{\frac{1}{r}})$.
\end{itemize}
\end{theorem}

\begin{proof}
The proof is similar to that of Theorem~\ref{glocal}, using Theorem~\ref{Spectral2} instead of Theorem~\ref{Spectral1}.
\end{proof}

The notion of hyperbolicity extends to any operator $T \in L(X)$ with exactly the same definition: 
$T \in L(X)$ is said to be {\em hyperbolic} if $\sigma(T) \cap \T = \emptyset$.
It is well known that $T \in L(X)$ is hyperbolic if and only if $X$ admits a direct sum decomposition $X = M \oplus N$, where 
$M$ and $N$ are closed $T$-invariant subspaces of $X$ such that $T|_N : N \to N$ is invertible, $r(T|_M) < 1$ and $r((T|_N)^{-1}) < 1$.
The concept of generalized hyperbolicity was extended to any $T \in L(X)$ in \cite[Proposition~19]{NBerAMes-Arxiv}, 
although the terminology ``generalized hyperbolic'' appeared only in \cite[Section~5]{BP}.
An operator $T \in L(X)$ is said to be {\em generalized hyperbolic} if $X$ admits a direct sum decomposition $X = M \oplus N$,
where $M$ and $N$ are closed subspaces of $X$ with the following properties:
\begin{itemize}
\item $T(M) \subset  M$ and $r(T|_M) < 1$;
\item $T|_N : N \to T(N)$ is bijective, $T(N)$ is closed, $N \subset T(N)$ and $r((T|_N)^{-1}|_N) < 1$.
\end{itemize}

\begin{corollary}[{\cite[Proposition 19]{NBerAMes-Arxiv}}]\label{GH-PSP}
Every generalized hyperbolic operator $T \in L(X)$ has the positive shadowing property.
\end{corollary}

\begin {proof}
Let $X = M \oplus N$ be a direct sum decomposition given by the generalized hyperbolicity of $T$.
Let $S = (T|_N)^{-1}|_N \in L(N)$ and choose a real number $r$ such that $\max\{r(T|_M),r(S)\} < r < 1$.
Clearly, $M \subset \{x \in X : r_T(x) \leq r\}$.
Hence, by Theorem~\ref{glocal2}, it remains to prove that
\begin{equation}\label{cor-glocal2-eq1}
N \subset \cX_T(\C \backslash \D_{\frac{1}{r}}).
\end{equation}
Fix $x \in N$ and define $f : \D_{\frac{1}{r}} \to N$ by $f(\lambda) = - (\lambda S - I_N)^{-1} Sx$,
where $I_N$ denotes the identity operator on $N$. 
Note that $f$ is well-defined, since $0 < |\lambda| < 1/r$ implies that $\lambda^{-1} \not\in \sigma(S)$
and, therefore, $\lambda S - I_N$ is invertible.
Moreover, $f$ is an analytic function on $\D_{\frac{1}{r}}$.
We claim that
\begin{equation}\label{cor-glocal2-eq2}
(T - \lambda I) f(\lambda) = x \ \ \text{ for all } \lambda \in \D_{\frac{1}{r}}.
\end{equation}
Indeed, for $\lambda = 0$, $(T - \lambda I) f(\lambda) = T f(0) = TSx = x$, and for $\lambda \in \D_{\frac{1}{r}} \backslash \{0\}$,
\begin{align*}
(T - \lambda I) f(\lambda) &= (T|_N - \lambda I_N)(- \lambda^{-1} (S - \lambda^{-1} I_N)^{-1} Sx)\\
  &= T|_N (S - \lambda^{-1} I_N) (S - \lambda^{-1} I_N)^{-1} Sx = x.
\end{align*}
By (\ref{cor-glocal2-eq2}), $x \in \cX_T(\C \backslash \D_{\frac{1}{r}})$, which proves (\ref{cor-glocal2-eq1}).
\end{proof}

The previous proof can be adapted to establish the following more general result.

\begin{corollary}
Let $T \in L(X)$. Suppose that there exists a mapping $B : X \to X$ (not necessarily continuous or linear) with the following properties:
\begin{itemize}
\item [\rm (c1)] $B(t x) = t B(x)$ for all $x \in X$ and all $t \geq 0$.
\item [\rm (c2)] $\displaystyle \limsup_{n \to \infty} \|T^n \circ B\|^{\frac{1}{n}} < 1$.
\item [\rm (c3)] If $N$ denotes the set $(I-B)(X)$, then $T|_N : N \to T(N)$ is bijective, $N \subset T(N)$ and the mapping
  $S = (T|_N)^{-1}|_N : N \to N$ satisfies $\displaystyle \limsup _{n \to \infty} \|S^n \circ (I-B)\|^{\frac{1}{n}} < 1$.
\end{itemize}
Then $T$ has the positive shadowing property.
\end{corollary}

Indeed, take $r < 1$ larger than the superior limits in (c2) and (c3). Clearly, $B(X) \subset \{x \in X : r_T(x) \leq r\}$.
Given $x \in N$, define $f : \D_{\frac{1}{r}} \to X$ by
\[
f(\lambda) = \sum_{n=0}^\infty \lambda^n S^{n+1}(x).
\]
Note that the series above converges and defines an analytic function on $\D_{\frac{1}{r}}$.
Elementary computations show that
\[
(T - \lambda I) f(\lambda) = x \ \ \text{ for all } \lambda \in \D_{\frac{1}{r}}.
\]
Hence, $x \in \cX_T(\C \backslash \D_{\frac{1}{r}})$. This proves that $(I-B)(X) = N \subset \cX_T(\C \backslash \D_{\frac{1}{r}})$.
Thus, by Theorem~\ref{glocal2}, $T$ has the positive shadowing property.

\medskip
An operator $T \in L(X)$ is said to have the {\em single-valued extension property} (SVEP) if, 
for any open subset $U$ of $\C$, the only analytic function $f: U \to X$ satisfying the equation 
$(T - \lambda I)f(\lambda) = 0$ identically on $U$ is the zero function on $U$.

An operator $T \in GL(X)$ is {\em expansive} if $\sup_{n \in \Z} \|T^nx\| = \infty$ for every nonzero $x \in  X$.

\begin{theorem}\label{thm-SVEP}
The following statements are true: 
\begin{itemize}
\item [\rm (a)] If $T \in L(X)$ has the single-valued extension property, then $T$ has the positive shadowing property
  if and only if $T$ is hyperbolic.
\item [\rm (b)] If $T \in GL(X)$ has the single-valued extension property, then $T$ has the shadowing property
  if and only if $T$ is hyperbolic.
\item [\rm (c)] If $T \in GL(X)$ is expansive, then $T$ has the shadowing property if and only if $T$ is hyperbolic.
\end{itemize}
\end{theorem}

\begin{proof}
(a): If $T \in L(X)$ has the single-valued extension property, 
then \cite [Proposition~1.3.2(f)]{LN} guarantees that $\sigma_s(T) =  \sigma(T)$. 
Therefore, by Theorem~\ref{Spectral2}, $T$ has the positive shadowing property if and only if $T$ is hyperbolic.

\smallskip\noindent
(b): It follows immediately from (a), since the concepts of shadowing and positive shadowing 
coincide for invertible operators on Banach spaces \cite[Theorem~1]{BP}.

\smallskip\noindent
(c): If $T \in GL(X)$ is expansive, then $\sigma_p(T) \cap \T = \emptyset$. 
If, in addition, $T$ has the shadowing property, then $\sigma_s(T) \cap \T = \emptyset$. 
Thus, $\sigma(T) \cap \T = \emptyset$, that is, $T$ is hyperbolic.
\end{proof}

An operator $T \in L(X)$ is said to be {\em paranormal} if $\|Tx\|^2 \leq \|T^{2}x\|\|x\|$ for all $x \in X$.

An operator $T \in L(X)$ is said to have {\em finite ascent} if, for every $\lambda \in \C$, 
there exists $n \in \N$ such that $\ker(T - \lambda I)^n = \ker(T - \lambda I)^{n+1}$.

\begin{corollary}
Let $T \in  L(X)$. If
\begin{itemize}
\item [\rm (a)] $\dim X < \infty$, or
\item [\rm (b)] $(\sigma_p(T))^{\rm o} =\emptyset$, or
\item [\rm (c)] $X$ is a Hilbert space and $T$ is paranormal, or
\item [\rm (d)] $X$ is a separable Banach space and $T$ is paranormal, or
\item [\rm (e)] $T$ has finite ascent, 
\end{itemize}
then T has the positive shadowing property if and only if T is hyperbolic.
\end{corollary}

\begin{proof}
In any of the cases above, the operator $T$ has the single-valued extension property.
Indeed, cases (a) and (b) are clear.
For cases (c), (d) and (e), see \cite[Page~2445]{AG}, \cite[Corollary~2.10]{CR} and \cite[Proposition~1.8]{L}, respectively.
\end{proof}

An operator $T \in L(X)$ is said to be {\em supercyclic} if there is a vector $x \in X$ whose projective orbit
$\{\lambda T^nx : \lambda \in \C, \ n \in \N_0\}$ is dense in $X$.

\begin{corollary}
If $T \in L(X)$ is a supercyclic operator, then $T^*$ has the positive shadowing property if and only if $T^*$ is hyperbolic.
\end{corollary}

\begin{proof}
If $T \in L(X)$ is a supercyclic operator, then $(\sigma_p(T^*))^{\rm o} = \emptyset$ \cite[Proposition~3.1]{H2}.
Therefore, the result follows from the previous corollary.
\end{proof}

Since the concepts of shadowing and hyperbolicity coincide for invertible operators in the finite-dimensional setting, 
it is natural to ask whether there exists an infinite-dimensional Banach space in which these concepts also coincide.
The following theorem provides an affirmative answer.

\begin{theorem}\label{SS}
There is an infinite-dimensional (separable) Banach space $X$ such that an operator $T \in L(X)$ (resp.\ $T \in GL(X)$)
has the positive shadowing property (resp.\ the shadowing property) if and only if it is hyperbolic.
\end{theorem}

\begin{proof}
Let $X$ be a complex hereditarily indecomposable (separable) Banach space. 
By \cite[Theorem~18]{GM}, every operator $T \in L(X)$ can be written as $T = \lambda I + S$, 
where $\lambda \in \C$ and $S$ is strictly singular;
moreover, the spectrum $\sigma (T)$ of $T$ is finite or consists of $\lambda$ 
and a sequence $(\lambda_n)_{n \in \N}$ of eigenvalues with finite multiplicity converging to $\lambda$. 
Thus, $\sigma(T) = \partial \sigma(T) = \sigma_s(T)$, 
and we conclude that $T$ has the positive shadowing property if and only if $T$ is hyperbolic.
\end{proof}

An operator $T \in L(X)$ is said to be {\em Li-Yorke chaotic} if there is an uncountable set $S \subset X$ such that
each pair $(x,y)$ of distinct elements of $S$ is a {\em Li-Yorke pair} for $T$, in the sense that
\[
\liminf_{n \to \infty} \|T^nx - T^ny\| = 0 \ \ \text{ and } \ \ \limsup_{n \to \infty} \|T^nx - T^ny\| > 0.
\]
$T$ is said to be {\em hypercyclic} if there is a vector $x \in X$ whose orbit $\{T^nx : n \in \N_0\}$ is dense in $X$.
It is well-known that every hypercyclic operator is Li-Yorke chaotic.
Since hyperbolic operators are never Li-Yorke chaotic, we obtain the following result.

\begin{corollary}
There exists an infinite-dimensional (separable) Banach space $X$ such that 
every operator $T \in L(X)$ with the positive shadowing property is not Li-Yorke chaotic (and, therefore, is not hypercyclic).
\end{corollary}

It was shown in \cite[Proposition 1]{CirGolPuj21} that if $T \in GL(X)$ is a  generalized hyperbolic operator, 
then either $T$ is hyperbolic or $\T \subset (\sigma_p(T))^{\rm o}$.
The next theorem shows that this result holds for operators with the shadowing property.

\begin{theorem}\label{PSP-H-or-PSpec}
If an operator $T\in L(X)$ has the positive shadowing property, then either $T$ is hyperbolic or $\T \subset (\sigma_p(T))^{\rm o}$.
In particular, the same conclusion holds for every $T \in GL(X)$ with the shadowing property.
\end{theorem}

\begin{proof}
Suppose that $T$ has the positive shadowing property and is not hyperbolic. Then,
\[
\sigma_s(T) \cap \T = \emptyset  \ \ \text{ and } \ \ \sigma(T) \cap \T \neq \emptyset.
\]
Note that we necessarily have $\T \subset \sigma(T)$, for otherwise there would exist $\lambda \in \T$ such that 
$\lambda \in \partial \sigma(T) \subset \sigma_s(T)$, which would constitute a contradiction. 
Therefore, 
\[
\T \subset \sigma_p(T),
\]
since $\sigma(T) = \sigma_p(T) \cup \sigma_s(T)$ and $\sigma_s(T) \cap \T = \emptyset $.
From this, we conclude that $\T \subset (\sigma_p(T))^{\rm o}$, since, otherwise, 
there would exist $\lambda \in \T$ with $\lambda \in \partial\sigma(T) \subset \sigma_s(T)$.
\end{proof}

An operator $T \in L(X)$ is said to have the {\em single-valued extension property} at $\lambda_0 \in \C$ (SVEP at $\lambda_0$) if,
for every neighborhood $U$ of $\lambda_0$, the only analytic function $f : U \to X$ satisfying the equation
$(T - \lambda I)f(\lambda) = 0$ identically on $U$ is the zero function on $U$.

\begin{corollary}
If $T \in L(X)$ is a non-hyperbolic operator with the positive shadowing property, then, for every $\lambda \in \T$,
$T$ does not have the single-valued extension property at $\lambda$.
\end{corollary}

\begin{proof}
If $T$ has the positive shadowing property, then $\lambda I - T$ is surjective for all $\lambda \in \T$. 
But, by \cite [Corollary~2.24]{A}, if $\lambda I - T$ is surjective, 
then $T$ has the SVEP at $\lambda$ if and only if $\lambda I - T$ is injective.
On the other hand, if $T$ is also non-hyperbolic, then the previous theorem guarantees that $\T \subset (\sigma_p(T))^{\rm o}$.
By combining these facts, we conclude that $T$ cannot have the SVEP at $\lambda$ if $\lambda \in \T$.
\end{proof}

Recall that $T \in L(X)$ is said to be {\em topologically transitive} (resp.\ {\em topologically mixing}) if, 
for any pair $U,V$ of nonempty open sets in $X$, there exists $n \in \N_0$ (resp.\ $n_0 \in \N_0$) such that 
$T^n(U) \cap V \neq \emptyset$ (resp.\ for all $n \geq n_0$).
If the space $X$ is separable, then Birkhoff's transitivity theorem asserts that $T$ is topologically transitive if and only if it is hypercyclic
(see \cite[Theorem~1.16]{GP}).
Recall also that $T$ is said to be {\em Devaney chaotic} if it is topologically transitive and the set of periodic points of $T$ is dense in $X$. 
Furthermore, $T$ is called {\em sub-chaotic} if there exists a closed $T$-invariant subspace $Y \neq \{0\}$ of $X$ 
such that $T|_Y$ is Devaney chaotic.

\begin{corollary}
If $T\in L(X)$ is a non-hyperbolic operator with the positive shadowing property, then $T$ is sub-chaotic.
\end{corollary}

\begin{proof}
If $T \in L(X)$ is a non-hyperbolic operator with the positive shadowing property, then
for any $\lambda \in \T$, we have that $\lambda \in \sigma_p(T)$ and $T - \lambda I$ is surjective.
Thus, by \cite[Theorem~2.2]{M}, $T$ is sub-chaotic.
\end{proof}

Recall that $T \in L(X)$ is said to be {\em distributionally chaotic} if there exist an uncountable set $S \subset X$ and $\eps > 0$
such that each pair $(x,y)$ of distinct elements of $S$ is an {\em $\eps$-distributionally chaotic pair} for $T$, in the sense that
\[
\udens\{n \in \N : \|T^n x - T^n y\| \geq \eps\} = 1, \text {and}
\]
\[
\udens\{n \in \N : \|T^n x - T^n y\| < \delta\} = 1 \ \text{ for all } \delta > 0,
\]
where $\udens(I)$ denotes the {\em upper density} of the subset $I$ of $\N$, that is,
\[
\udens(I) = \limsup_{n \to \infty} \frac{\card(I \cap [1,n])}{n}\cdot
\]
If $S$ can be chosen to be dense in $X$, then $T$ is {\em densely distributionally chaotic}.
Finally, by a {\em distributionally irregular manifold} for $T$, we mean a vector subspace $Y$ of $X$ such that 
every non-zero vector $y \in Y$ is a {\em distributionally irregular vector} for $T$, that is, 
there exist sets $A,B \subset \N$ with $\udens(A) = \udens(B) = 1$ such that
$\lim_{n \in A} \|T^n y\| = 0$ and $\lim_{n \in B} \|T^n y\| = \infty$.

\begin{corollary}
If $T\in L(X)$ is a non-hyperbolic operator with the positive shadowing property and there exists a dense subset $X_0$ of $X$ such that 
\[
T^nx \to 0 \ \ \text{ for all } x \in X_0,
\]
then $T$ is distributionally chaotic. In, in addition, $X$ is separable, then $T$ has a dense distributionally irregular manifold
and, in particular, is densely distributionally chaotic.
\end{corollary}

\begin{proof}
Since $T$ is non-hyperbolic and has the positive shadowing property, Theorem~\ref{PSP-H-or-PSpec} shows that 
there exists $\lambda \in \sigma_p(T)$ with $|\lambda| > 1$. 
Hence, the result follows from \cite[Theorem~33]{BBPW}.
\end{proof}

Recall that $T \in L(X)$ is said to be {\em mean Li-Yorke chaotic} if there exists an uncountable set $S \subset X$ such that 
each pair $(x,y)$ of distinct points in $S$ is a {\em mean Li-Yorke pair} for $T$, that is,
\[
\liminf_{N \to \infty} \frac{1}{N} \sum_{j=1}^N \|T^j x - T^j y\| = 0
\ \ \text {and} \ \
\limsup_{N \to \infty} \frac{1}{N} \sum_{j=1}^N \|T^j x - T^j y\| > 0.
\]
If $S$ can be chosen to be dense in $X$, then $T$ is {\em densely mean Li-Yorke chaotic}.
Finally, by an {\em absolutely mean irregular manifold} for $T$, we mean a vector subspace $Y$ of $X$ such that 
every non-zero vector $y \in Y$ is an {\em absolutely mean irregular vector} for $T$, that is, 
\[
\liminf_{N \to \infty} \frac{1}{N} \sum_{j=1}^N \|T^j y\| = 0
\ \ \text{ and } \ \
\limsup_{N \to \infty} \frac{1}{N} \sum_{j=1}^N \|T^j y\| = \infty.
\]

\begin{corollary}
If $T\in L(X)$ is a non-hyperbolic operator with the positive shadowing property and there exists a dense subset $X_0$ of $X$ such that 
\[
\lim_{N \to \infty} \frac{1}{N} \sum_{j=1}^N \|T^jx\| = 0 \ \ \text{ for all } x \in X_0,
\]
then $T$ is mean Li-Yorke chaotic. If, in addition, $X$ is separable, then $T$ has a dense absolutely mean irregular manifold
and, in particular, is densely mean Li-Yorke chaotic.
\end{corollary}

\begin{proof}
Since $T$ is hyperbolic and has the positive shadowing property, it follows from Theorem~\ref{PSP-H-or-PSpec} that 
there exists $\lambda \in \sigma_p(T)$ with $|\lambda| > 1$. 
Therefore, the result follows from \cite[Theorem~22]{BBPe} and \cite[Theorem~29]{BBPe}.
\end{proof}

Given $T \in L(X)$, $x,y \in X$ and $\delta > 0$, recall that a {\em $\delta$-chain for $T$ from $x$ to $y$}
is a finite $\delta$-pseudotrajectory $(x_n)_{n=0}^k$ of $T$ with $x_0 = x$ and $x_k = y$.
The operator $T$ is said to be {\em chain recurrent} (resp.\ {\em chain transitive}) if, 
for every $x \in X$ (resp.\ $x,y \in X$) and every $\delta > 0$, there is a $\delta$-chain for $T$ from $x$ to itself (resp.\ from $x$ to $y$).

It is known that if an operator $T$ on a separable Banach space $X$ is chain recurrent and has the positive shadowing property, 
then $T$ is topologically mixing, frequent hypercyclic and densely distributionally chaotic (see \cite{AMV} and \cite{BP}).
As a consequence of the next theorem, we shall prove that such an operator is also densely mean Li-Yorke chaotic.

\begin{theorem}
If an operator $T$ on a Banach space $X$ has the positive shadowing property, 
then $T$ has a residual set of absolutely mean irregular vectors if and only if the following conditions hold:
\begin{itemize}
\item [\rm (a)] The set $\{x \in X : \text{for each } \delta > 0, \text{ there is a $\delta$-chain for $T$ from $x$ to $0$}\}$ is dense in $X$.
\item [\rm (b)] There exists a vector $y \in X$ such that $\,\displaystyle \limsup_{N \to \infty} \frac{1}{N} \sum_{j=1}^N \|T^jy\| > 0$.
\end{itemize}
If, in addition, $X$ is separable, then $T$ is densely mean Li-Yorke chaotic.
\end{theorem}

\begin{proof}
Let $A$ be the set of all absolutely mean irregular vectors for $T$ and let $B$ denote the set described in condition (a).

Suppose that conditions (a) and (b) are true.
We claim that, for each $t \in \N$, the set
\[
B_t = \bigcup_{N=1}^\infty \Big\{x \in X : \frac{1}{N} \sum_{j=1}^N \|T^j x\| < \frac{1}{t} \Big\}
\]
is dense in $X$. Indeed, take $x \in B$ and $\eps \in (0,\frac{1}{2t})$.
Let $\delta > 0$ be associated to $\eps$ according to the positive shadowing property of $T$.
Since $x \in B$, there is a $\delta$-chain $(x_j)_{j=0}^k$ for $T$ from $x$ to $0$.
Since the sequence
\[
(x_j)_{j \in \N_0} = (x_0,x_1,\ldots,x_k,0,0,\ldots)
\]
is a $\delta$-pseudotrajectory of $T$, there exists $y \in X$ such that
\[
\|x_j - T^j y\| < \eps \ \ \text{ for all } j \in \N_0.
\]
Hence, $\|x - y\| < \eps$ and $y \in B_t$ since, for every sufficiently large $N \in \N$,
\[
\frac{1}{N} \sum_{j=1}^N \|T^j y\| < \frac{\|Ty\| + \cdots + \|T^ky\|}{N} + \frac{N-k}{N}\, \eps < 2\eps < \frac{1}{t}\cdot
\]
This proves our claim, since $B$ is dense in $X$. As a consequence, the set
\[
\Big\{x \in X : \liminf_{N \to \infty} \frac{1}{N} \sum_{j=1}^N \|T^j x\| = 0 \Big\} \ \text{ is residual in } X.
\]
In view of \cite[Theorem~22]{BBPe}, this fact, together with condition (b), implies that the set $A$ is residual in~$X$.

Conversely, suppose that $A$ is residual in $X$.
If $x \in A$, then
\[
\liminf_{N \to \infty} \frac{1}{N} \sum_{j=1}^N \|T^n x\| = 0,
\]
which implies that $\inf_{n \in \N} \|T^n x\| = 0$ and, therefore, $x \in B$.
Hence, $A \subset B$ and condition (a) holds.
On the other hand, every $y \in A$ satisfies the property described in (b).
Thus, condition (b) also holds.

Finally, the last assertion in the theorem follows from \cite[Theorem~17]{BBPe}.
\end{proof}

\begin{corollary}
If an operator $T$ on a Banach space $X$ is chain recurrent and has the positive shadowing property,
then $T$ has a residual set of absolutely mean irregular vectors.
If, in addition, $X$ is separable, then $T$ is densely mean Li-Yorke chaotic.
\end{corollary}

\begin{proof}
Since $T$ is chain recurrent and the concepts of chain recurrence and chain transitivity coincide for linear operators
\cite[Proposition~19]{ABM}, condition~(a) in the previous theorem holds.
It remains to establish condition~(b).
For this purpose, let $\delta > 0$ be associated to the number $1$ according to the fact that $T$ has the positive shadowing property.
Choose a vector $x \in X$ with $\|x\| > 1$.
Since $T$ is chain recurrent, there is a $\delta$-chain $(x_j)_{j=0}^k$ for $T$ from $x$ to itself.
Since
\[
(x_0,x_1,\ldots,x_k,x_1,\ldots,x_k,x_1,\ldots,x_k,\ldots)
\]
is a $\delta$-pseudotrajectory of $T$, there exists $y \in X$ such that
\[
\|x - T^{nk} y\| < 1 \ \ \text{ for all } n \in \N_0.
\]
Thus,
\[
\limsup_{N \to \infty} \frac{1}{N} \sum_{j=1}^N \|T^j y\| \geq \frac{\|x\| - 1}{k} > 0,
\]
as desired.
\end{proof}

Recall that $T \in L(X)$ is said to be {\em recurrent} if, for every nonempty open set $U$ in $X$,
there exists $n \in \N$ such that $U \cap T^{-n}(U) \neq \emptyset$.

\begin{theorem}
There exists an invertible non-hyperbolic operator with the shadowing property that is not recurrent.
\end{theorem}

\begin{proof}
Consider $0 < r < s < 1$. Take an invertible hyperbolic operator $S$ with $\sigma(S) \subset \D_r$ 
and a non-hyperbolic operator $R$ with the shadowing property such that $\sigma(R) \subset \C \backslash \D_s$. 
Then $T = S \oplus R$ verifies the result, since $\sigma(T) \cap \T \neq \emptyset$, $\sigma_s(T) \cap \T = \emptyset$,
and there exists a connected component of $\sigma(T)$ that does not intersect $\T$, 
which implies that $T$ is not recurrent \cite[Proposition~2.11]{CMP}.
\end{proof}

Given $T \in L(X)$, the {\em modulus of injectivity} of $T$ is defined by 
\[
j(T) = \inf\{\|Tx\| : \|x\| = 1\}.
\]
Let 
\[
i(T) = \lim_{n \to \infty} j(T^n)^{\frac{1}{n}}.
\]
Note that $i(T) \leq r(T)$.

\begin{theorem} 
Let $T \in L(X)$ be a non-invertible operator. 
\begin{itemize}
\item [\rm (a)] If $i(T) > 1$, then $T^*$ has the positive shadowing property and is not hyperbolic.
\item [\rm (b)] If $i(T) = r(T)$, then $T$ has the positive shadowing property if and only if $T$ is hyperbolic if and only if $r(T) < 1$.
\item [\rm (c)] If $i(T) = r(T)$, then $T^*$ has the positive shadowing property if and only if $r(T) \neq 1$.
\end{itemize}
\end{theorem}

\begin{proof}
(a): Since $T$ is a non-invertible operator, \cite[Theorem~2.55]{A} tells us that
\[
\ov{\D}_{i(T)} \subset \sigma(T). 
\]
Hence, $i(T) > 1$ implies that $\sigma(T^*) \cap \T \neq \emptyset$ and, therefore, $T^*$ is not hyperbolic.

On the other hand, by \cite[Theorem~2.54]{A},
\[
\sigma_s(T^*) \subset \{\lambda \in \C : i(T) \leq |\lambda| \leq r(T)\}.
\]
If $i(T) > 1$, then $\sigma_s(T^*) \cap \T = \emptyset$ and, therefore, $T^*$ has the positive shadowing property.

\smallskip\noindent
(b,c): Since $T$ is a non-invertible operator and $i(T) = r(T )$, it follows from \cite[Theorem 4.68]{A2} that 
\[
\sigma(T) = \ov{\D}_{r(T)} \ \ \text{ and } \ \ \sigma_s(T^*) = \partial \D_{r(T)}.
\]
The second equality clearly implies the equivalence in (c) and the first inequality clearly implies the second equivalence in (b).
Both equalities together imply that $\sigma_s(T) = \sigma(T)$, which gives the first equivalence in (b).
\end{proof}

\begin {remark}
For every $T \in GL(X)$, \cite[Theorem~2.55]{A} tells us that
\[
\sigma (T) \subset \{\lambda \in \C : i(T) \leq |\lambda| \leq r(T)\}. 
\]
Thus, if $i(T) > 1$ or $r(T) < 1$, then $\sigma(T) \cap \T = \emptyset$ and $T$ is  hyperbolic.
\end{remark}

An operator $T \in GL(X)$ is said to be {\em uniformly expansive} if there exists $n \in \N$ such that,
for each $x \in X$ with $\|x\| = 1$, we have that $\|T^nx\| \geq 2$ or $\|T^{-n}x\| \geq 2$.

Let $H$ be a complex Hilbert space.
For each operator $T \in L(H)$, the associated operators $L_T, R_T \in L(L(H))$ are defined by $L_T(A) = T A$ and $R_T(A) = A T$.

\begin{theorem} 
Let $T \in GL(H)$, where $H$ is a complex Hilbert space.
\begin{itemize}
\item [\rm (a)] If $T$ has the shadowing property on $H$, then $L_T$ has the shadowing property on $L(H)$.
\item [\rm (b)] If $T$ is uniformly expansive on $H$, then $R_T$ has the shadowing property on $L(H)$.
\end{itemize}
\end{theorem}

\begin{proof}
(a): Since $\sigma_s(L_T) \subset \sigma_s(T)$ \cite[Page~257]{LN}, the result follows from Theorem~\ref{Spectral1}.

\smallskip\noindent
(b): Since $\sigma_s(R_T) \subset \sigma_a(T)$ \cite[Page~257]{LN}, the result follows from Theorem~\ref{Spectral1}
and the fact that $T$ is uniformly expansive if and only if $\sigma_a(T) \cap \T = \emptyset$ \cite[Theorem~2]{EH}.
\end{proof}

Given a {\em domain} (nonempty connected open set) $\Omega \subset \C$ and a positive integer $n$, 
Cowen and Douglas \cite{CD} introduced the operator class $\cB_n(\Omega)$, 
consisting of all operators $T$ on a complex separable Hilbert space $H$ that satisfy:

1. $\Omega \subset \sigma(T)$;

2. $\ran(T - \omega I) = H$ for $\omega \in \Omega$;

3. $\spa _{\omega \in \Omega} \ker(T - \omega I) = H$; and

4. $\dim \ker(T - \omega I) = n$ for $\omega \in \Omega$.

\begin{theorem}
Given a domain $\Omega \subset \C$, a positive integer $n$ and an operator $T \in \cB_n(\Omega)$, we have that:
\begin{itemize}
\item [\rm (a)] If $\T \subset \Omega$, then $T$ is a non-hyperbolic operator with the positive shadowing property.
\item [\rm (b)] $T^*$ has the positive shadowing property if and only if it is hyperbolic.
\end{itemize}
\end{theorem}

\begin{proof}
(a): If $T \in \cB_n(\Omega)$ and $\T \subset \Omega$, then $\sigma(T) \cap \T \neq \emptyset$ and $\sigma_s(T) \cap \T = \emptyset$.

\smallskip\noindent
(b): If $T \in \cB_n(\Omega)$, then $\sigma_{p}(T^*) = \emptyset$ \cite[Theorem~1.1]{H}.
Thus, $T^*$ has the positive  shadowing property if and only if it is hyperbolic.
\end{proof}

Let $H$ be a complex separable infinite-dimensional Hilbert space and let $U \in  L(H)$. 
We say that $U$ is {\em universal in the sense of Caradus} \cite{C} if 

1. $\ker(U)$ is infinite-dimensional; and

2. $\ran(U)$ is $H$.

\begin{theorem}
If $T \in L(H)$ and, for each $\lambda \in \T$, the operator $T - \lambda I$ is universal in the sense of Caradus,
then $T$ is a non-hyperbolic operator with the positive shadowing property.
\end{theorem}

\begin{proof}
If, for each $\lambda \in \T$, the operator $T - \lambda I$ is universal in the sense of Caradus, 
then $\sigma(T) \cap \T \neq \emptyset$ and $\sigma_s(T) \cap \T = \emptyset$.
\end{proof}


\section {Operators  without the shadowing property}\label{Section4}

We now summarize a list of types of operators that never have the positive shadowing property.

\begin{theorem}
The following classes of operators on Banach spaces do not contain any operator with the positive shadowing property:
\begin{itemize}
\item [\rm (a)] Operators $T$ with $\partial \sigma(T) \cap \T \neq \emptyset$.
\item [\rm (b)] Operators $T$ with $r(T) = 1$.
\item [\rm (c)] $m$-isometries.
\item [\rm (d)] Operators $T$ with $\sigma_a(T^*) \cap \T \neq \emptyset$.
\item [\rm (e)] Li-Yorke chaotic operators with $\sigma(T) \subset \ov{\D}$.
\item [\rm (f)] Recurrent operators with $\sigma(T) \subset \ov{\D}$.
\item [\rm (g)] Adjoints of recurrent operators.
\item [\rm (h)] Hypercyclic compact perturbations of a multiple of the identity.
\end{itemize}
\end{theorem}

\begin{proof}
(a): It follows from the fact that $\partial \sigma(T) \subset \sigma_s(T)$.

\smallskip\noindent
(b): It follows from (a), since $r(T) = 1$ implies that $\partial \sigma(T) \cap \T \neq \emptyset$.

\smallskip\noindent
(c): It follows from (b), since every $m$-isometry $T$ satisfies $r(T) = 1$.

\smallskip\noindent
(d): It follows from the fact that $\sigma_a(T^*) = \sigma_s(T)$.

\smallskip\noindent
(e,f): If $T$ is Li-Yorke chaotic or recurrent, then $\sigma(T) \cap \T \neq \emptyset$.
If, in addition, $\sigma(T) \subset \ov{\D}$, then $r(T) = 1$ and we fall back into case (b).

\smallskip\noindent
(g): Let $T$ be a recurrent operator. 
On the one hand, \cite[Proposition~2.14]{CMP} guarantees that $\sigma_p(T^*) \subset \T$ and, therefore,
\[
(\sigma_p(T^*))^{\text{o}} = \emptyset.
\]
On the other hand, we have already mentioned that $\sigma(T) \cap \T \neq \emptyset$, 
that is, $\sigma(T^*) \cap \T \neq \emptyset$. Hence, 
\[
T^* \text{ is not hyperbolic}.
\]
By Theorem~\ref{PSP-H-or-PSpec}, we conclude that $T^*$ cannot have the positive shadowing property.

\smallskip\noindent
(h): Let $T$ be a hypercyclic operator of the form $T = \lambda I + K$, where $\lambda \in \C$ and $K$ is compact.
By \cite[Lemma~5.19]{GP}, $|\lambda| = 1$ and $\sigma(T) = \{\lambda\}$. 
Hence, we can apply case (a).
\end{proof}

\begin{theorem}
There is an infinite-dimensional (separable) Banach space $X$ such that no Li-Yorke chaotic operator 
(in particular, no hypercyclic operator) $T \in L(X)$ has the positive shadowing property.
\end{theorem}

\begin{proof}
Take the space $X$ as in Theorem~\ref{SS} and use the fact that Li-Yorke chaotic operators cannot be hyperbolic.
\end{proof}

Let us now mention some examples of operators that do not have the positive shadowing property.

\begin{example}
{\bf (a)} Let $X = C[0,1]$ of $X = L^p[0,1]$ for some $p \in [1,\infty)$.
The operator $I + V$, where $V$ is the {\em Volterra operator} on $X$ (see \cite[Example~5.4]{GP}), 
does not have the positive shadowing property, since $\sigma(V) = \{0\}$.

\smallskip\noindent
{\bf (b)} For any $p \in [1,\infty)$, if $B$ is the backward shift on the Banach space $\ell^p(\N)$ 
and $T = 2B \oplus I$ on $\ell^p(\N) \oplus \C$, then $T$ is a supercyclic operator that does not have the positive shadowing property.

\smallskip\noindent
{\bf (c)} For each $p \in [1,\infty)$ and each $\eps \in (0,1)$, 
it was shown in \cite[Corollary~2.3]{BBMP} that there exists an absolutely C\`esaro bounded operator $T$ on $\ell^p(\N)$ 
which is topologically mixing and satisfies
\begin{equation}\label{Example-Eq1}
\|T^n\| = (n+1)^{\frac{1-\eps}{p}} \ \ \text{ for all } n \in \N.
\end{equation}
By (\ref{Example-Eq1}), $r(T) = 1$. Therefore, $T$ does not have the positive shadowing property.

\smallskip\noindent
{\bf (d)} Given any non-empty closed subset $E$ of $\T$, it was shown in \cite[Corollary~4.8]{BG} that 
there exists a frequently hypercyclic operator $T$ on a separable Hilbert space such that 
\[
\sigma_p(T) = \sigma(T) = \sigma_s(T) = E.
\] 
Thus, $T$ does not have the positive shadowing property.
\end{example}


\section{Shadowing implies average shadowing}\label{Section5}

Since the introduction of the concept of average shadowing by Blank~\cite{Blank1}, 
various variations of this concept---involving average error in tracing---have been presented and investigated by several authors. 
We refer the reader to \cite{ABM,Blank2,WuOprChe16}, where many additional references can be found.

In our context, take $T \in GL(X)$.
A sequence $(x_n)_{n \in \Z}$ in $X$ is said to be:
\begin{itemize}
\item a {\em $\delta$-average-pseudotrajectory} of $T$, where $\delta > 0$, if there exists $N \in \N$ such that
  \[
  \frac{1}{n} \sum_{j=0}^{n-1} \|x_{k+j+1} - T x_{k+j}\| \leq \delta \ \text{ whenever } n \geq N \text{ and } k \in \Z;
  \]
\item an {\em asymptotic average pseudotrajectory} of $T$ if
  \[
  \lim_{n \to\infty} \Big(\frac{1}{2n+1} \sum_{j=-n}^{n} \|x_{j+1} - T x_j\|\Big) = 0.
  \]
\end{itemize}
The operator $T$ is said to have:
\begin{itemize}
\item the {\em average shadowing property} if, for every $\eps > 0$, there exists $\delta > 0$ such that 
  every $\delta$-average-pseudotrajectory $(x_n)_{n \in \Z}$ of $T$ is {\em $\eps$-shadowed on average} by some $x \in X$, that is,
  \[
  \limsup_{n \to\infty} \Big(\frac{1}{2n+1} \sum_{j=-n}^{n} \|x_j - T^j x\|\Big) < \eps;
  \]
\item the {\em asymptotic average shadowing property} if every asymptotic average pseudotrajectory $(x_n)_{n \in \Z}$ of $T$ 
  is {\em asymptotically shadowed on average} by some $x \in X$, that is,
  \[
  \lim_{n \to\infty} \Big(\frac{1}{2n+1} \sum_{j=-n}^{n} \|x_j - T^j x\|\Big) = 0.
  \]
\end{itemize}

In the setting of linear dynamics, the average shadowing property was studied in the articles \cite{ABM,JiaWanLi19}.
It was shown in \cite[Theorem~3]{JiaWanLi19} (see also \cite[Proposition~13]{ABM}) that every generalized hyperbolic operator 
$T \in GL(X)$ has the average shadowing property and the asymptotic average shadowing property.
The next theorem generalizes this result. 
The main strategy in its proof is an adaptation of the original proof, presented in \cite[Theorem~A]{BCDMP}, 
that generalized hyperbolicity implies the shadowing property.

\begin{theorem}
Let $T\in  GL(X)$. If there exists a mapping $B : X \to X$ (not necessarily continuous or linear) such that
$B(t x) = t B(x)$ for all $x \in X$ and all $t \geq 0$,
\[
\limsup_{n \to \infty} \|T^n \circ B\|^{\frac{1}{n}} < 1 \ \ \text{ and } \ \ \limsup _{n \to \infty} \|T^{-n} \circ (I-B)\|^{\frac{1}{n}} < 1,
\]
then $T$ has the average shadowing property and the asymptotic average shadowing property.
\end{theorem}

\begin{proof}
The hypotheses on $B$ imply the existence of constants $c \geq 1$ and $r \in (0,1)$ such that
\[
\|T^n(B(x))\| \leq c\, r^n \|x\| \ \text{ and } \ \|T^{-n}((I-B)(x))\| \leq c\, r^n \|x\| \ \ \ (x \in X, n \in \N_0).
\]

Let $(x_n)_{n \in \Z}$ be a sequence in $X$.
For each $n \in \Z$, let $z_n = x_{n+1} - Tx_n$.
Assume that the sequence $(z_n)_{n \in \Z}$ is at most $O(|n|)$, that is,
there is a constant $M \in (0,\infty)$ such that $\|z_n\| \leq M |n|$ for all $n \in \Z \backslash \{0\}$.
Then, for each $n \in \Z$, the vector
\[
y_n = \sum_{k=0}^\infty T^k(B(z_{n-k-1})) - \sum_{k=1}^\infty T^{-k}((I-B)(z_{n+k-1}))
\]
is well-defined, since the series above are absolutely convergent. 
Furthermore,
\begin{equation}\label{AverageEq1}
\|y_n\| \leq \sum_{k=0}^\infty c\, r^k \|z_{n-k-1}\| + \sum_{k=1}^\infty c\, r^k \|z_{n+k-1}\|.
\end{equation}
Since
\begin{align*}
z_n + T y_n &= z_ n + \sum_{k=0}^\infty T^{k+1}(B(z_{n-k-1})) - \sum_{k=1}^\infty T^{-k+1}((I-B)(z_{n+k-1}))\\
  &= z_n + \sum_{k=1}^\infty T^k(B(z_{n-k})) - \sum_{k=0}^\infty T^{-k}((I-B)(z_{n+k}))\\
  &= \sum_{k=0}^\infty T^k(B(z_{n-k})) - \sum_{k=1}^\infty T^{-k}((I-B)(z_{n+k}))\\
  &= y_{n+1},
\end{align*}
we obtain $x_{n+1} - y_{n+1} = T(x_n - y_n)$ for all $n \in \Z$.
Thus, $x_n - y_n = T^n(x_0 - y_0)$ and, therefore,
\begin{equation}\label{AverageEq2}
\|x_n - T^n(x_0 - y_0)\| = \|y_n\| \ \ \text{ for all } n \in \Z.
\end{equation}

To prove that $T$ has the average shadowing property, take $\eps > 0$, define $\delta = \frac{(1-r)\eps}{4c}$,
and assume that $(x_n)_{n \in \Z}$ is a $\delta$-average-pseudotrajectory of $T$.
In this case, $(z_n)_{n \in \Z}$ is a bounded sequence.
In view of (\ref{AverageEq2}), it remains to show that
\[
\limsup_{n \to\infty} \Big(\frac{1}{2n+1} \sum_{j=-n}^{n} \|y_j\|\Big) < \eps.
\]
Let $N \in \N$ be such that
\[
\frac{1}{2n+1} \sum_{j=-n}^{n} \|z_{j+k}\| \leq \delta \ \text{ whenever } n \geq N \text{ and } k \in \Z.
\]
By (\ref{AverageEq1}),
\begin{align}\label{AverageEq3}
\frac{1}{2n+1} \sum_{j=-n}^{n} \|y_j\| &\leq c \sum_{k=0}^\infty \Big(\frac{1}{2n+1} \sum_{j=-n}^{n} \|z_{j-k-1}\|\Big) r^k 
  + c \sum_{k=1}^\infty \Big(\frac{1}{2n+1} \sum_{j=-n}^{n} \|z_{j+k-1}\|\Big) r^k\\
  &\leq c \sum_{k=0}^\infty \delta\, r^k + c \sum_{k=1}^\infty \delta\, r^k < \frac{2\, c\, \delta}{1 - r} = \frac{\eps}{2} < \eps,\notag
\end{align}
for all $n \geq N$.

To prove that $T$ has the asymptotic average shadowing property, 
assume that $(x_n)_{n \in \Z}$ is an asymptotic average pseudotrajectory of $T$.
In this case, the sequence $(z_n)_{n \in \Z}$ is at most $O(|n|)$.
By (\ref{AverageEq2}), it suffices to prove that
\[
\lim_{n \to\infty} \Big(\frac{1}{2n+1} \sum_{j=-n}^{n} \|y_j\|\Big) = 0.
\]
For this purpose, take $\eps > 0$ and define
\[
s = \sum_{k=1}^\infty k\, r^k \ \ \ \text{ and } \ \ \ t_n = \frac{1}{2n+1} \sum_{j=-n}^n \|z_j\| \ \ \ (n \in \N).
\]
Let $n_0 \geq 2$ be such that $t_n < \eps$ for all $n \geq n_0$.
For every $n \geq n_0$,
\begin{align}\label{AverageEq4}
\sum_{k=0}^\infty \Big(\frac{1}{2n+1} \sum_{j=-n}^{n} \|z_{j-k-1}\|\Big) r^k
  &\leq \sum_{k=0}^\infty \Big(\frac{2n+2k+3}{2n+1}\, t_{n+k+1}\Big) r^k \notag\\
  &\leq 2\, t_{n+1} + \sum_{k=1}^\infty 2\, k\, t_{n+k+1}\, r^k < 2 (1 + s) \eps.
\end{align}
Similarly, for every $n \geq n_0$,
\begin{align}\label{AverageEq5}
\sum_{k=1}^\infty \Big(\frac{1}{2n+1} \sum_{j=-n}^{n} \|z_{j+k-1}\|\Big) r^k
  &\leq \sum_{k=1}^\infty \Big(\frac{2n+2k-1}{2n+1}\, t_{n+k-1}\Big) r^k \notag\\
  &\leq \sum_{k=1}^\infty k\, t_{n+k-1}\, r^k <  s\, \eps.
\end{align}
By (\ref{AverageEq3}), (\ref{AverageEq4}) and (\ref{AverageEq5}),
\[
\frac{1}{2n+1} \sum_{j=-n}^{n} \|y_j\| < c\, (2 + 3 s)\, \eps \ \ \text{ for all } n \geq n_0.
\]
Since $\eps > 0$ is arbitrary, we obtain the desired result. 
\end{proof}

\begin{corollary}
Every operator $T \in GL(X)$ with the shadowing property has the average shadowing property and the asymptotic average shadowing property.
\end{corollary}

\begin{proof}
It follows trivially from the previous theorem that every generalized pseudo-hyperbolic operator $T \in GL(X)$
has the average shadowing property and the asymptotic average shadowing property.
Since the shadowing property is equivalent to generalized pseudo-hyperbolicity \cite[Theorem~3.2]{DP2},
we obtain the desired result.
\end{proof}


\section{Shadowing for invertible weighted composition operators on the unit disk}\label{Section6}

In this section, we will study the shadowing property for invertible {\em weighted composition operators}
\[
(uC_\varphi f)(z) = u(z)(f \circ \varphi)(z)
\]
on certain Banach spaces $X$ of holomorphic functions on the unit disk $\D$, where $\varphi$ is an automorphism of $\D$
and $u$ is a multiplier of $X$.
Two special cases deserve mention: if $u$ is the constant function $1$, we obtain the {\em composition operator} $C_\varphi$,
and if $\varphi$ is the identity automorphism of $\D$, we obtain the {\em multiplication operator} $M_u$.

Recall that an {\em automorphism} of $\D$ is a bijective holomorphic map $\varphi : \D \to \D$.
We denote by $Aut(\D)$ the group of all automorphisms of $\D$.
It is well known that the non-trivial automorphisms $\varphi$ of $\D$ (i.e., $\varphi$ is not the identity automorphism $I_{\D}$)
fall into three disjoint classes, which can be characterized by the configuration of fixed points:
\begin{itemize}
\item $\varphi$ is {\em elliptic} if it has a unique fixed point in $\D$;
\item $\varphi$ is {\em parabolic} if it has a unique fixed point in $\T$;
\item $\varphi$ is {\em hyperbolic} if it has two distinct fixed points in $\T$.
\end{itemize}
Concerning the boundary fixed points of $\varphi$, let us recall the celebrated Denjoy–Wolff Theorem \cite[Theorem~2.51]{CM}:
{\it If $\varphi$ is either a parabolic or a hyperbolic automorphism of $\D$, then there is a (unique) fixed point $a \in \T$ such that
$\lim_{n \to \infty} \varphi^n(z) = a$ uniformly on compact subsets of $\D$ 
(the point a is called the {\em Denjoy–Wolff point} of $\varphi$).
Moreover, if $\varphi$ is parabolic, then $\varphi'(a) = 1$; 
and if $\varphi$ is hyperbolic and its other fixed point is $b \in \T$, then $0 < \varphi'(a) < 1$ and $\varphi'(b) = \frac{1}{\varphi'(a)}\cdot$}


\subsection{Shadowing for invertible weighted composition operators with hyperbolic symbol}

For any open subset $U$ of $\C$, we denote by $H(U)$ the algebra of all holomorphic functions on $U$.
Given a Banach space $X$ with $X \subset H(\D)$, let $Mul(X)$ be the algebra of all {\em multipliers} of $X$, 
that is, a function $g$ belongs to $Mul(X)$ when the pointwise product operator $f \to gf $ defines a bounded operator on $X$. 

The following definition was given in \cite[Definition~2.1]{Oliva}.

\begin{definition} 
Let $\gamma \geq 0$ and let $X \subset H(\D)$ be a non-zero Banach space of holomorphic functions. 
We say that $X$ is a {\em $\gamma$-space} if the pair $(X,\gamma)$ satisfies the following properties:
\begin{itemize}
\item [\rm (P1)] For each open set $U$ in $\C$ with $\ov{\D} \subset U$, we have that $H(U) \subset Mul(X)$.

\item [\rm (P2)] For each $c \in \T$ and each $\delta \in \C$ with $\re \delta > 0$, the function $z \to (c-z)^{\delta}$ belongs to $Mul(X)$.

\item [\rm (P3)] Let $\varphi$ be a hyperbolic automorphism with fixed points $a$ (attractive) and $b$ (repulsive) in $\T$. 
  For each $u \in Mul(X)$ for which $|u|$ has a continuous extension to $\D \cup \{a,b\}$, we have that
  the function $u_ n = \prod_{j=0}^{n-1} (u \circ \varphi^j)$ belongs to $Mul(X)$ for every $n \in \N$, and
  \[
  \limsup_{n \to \infty} \|u_n\|^{1/n}_{Mul(X)} \leq \max\{|u|(a),|u|(b)\},
  \]
  where $\|u_n\|_{Mul(X)}$ denotes the norm of the multiplier $u_n$ regarded as an operator on $X$.

\item [\rm (P4)] $(\varphi')^{\gamma} C_{\varphi} \in L(X)$ for all $\varphi \in Aut(\D)$, and for each $\eps > 0$, 
\[
\sup_{\varphi \in Aut(\D)} (1-|\varphi(0)|^2)^{\eps} \|(\varphi')^{\gamma} C_{\varphi}\|_{L(X)} < \infty.
\]
\end{itemize}
\end{definition}

As main examples of $\gamma$-spaces, we have: 
the Hardy spaces $H^p(\D)$, $1 \leq p \leq \infty$, with $\gamma = 1/ p$; 
the weighted Bergman spaces $A^p_{\sigma}(\D)$, $1 \leq p < \infty$, $\sigma > -1$, with $\gamma = \frac{\sigma+2}{p}$;  
the weighted Dirichlet spaces $D^p_{\sigma}(\D)$, $1 \leq p < \infty$, $\sigma > -1$, $\sigma \geq p-2$, with $\gamma = \frac{\sigma+2}{p}-1$; 
the disk algebra $A(\D)$, the Bloch space $\cal B$ and little Bloch space ${\cal B}_0$ with $\gamma = 0$.

\begin{theorem} 
Let $X$ be a $\gamma$-space for some $\gamma \geq 0$ and let $uC_{\varphi}$ be an invertible weighted composition operator on $X$.
Assume that the automorphism $\varphi$ is hyperbolic with fixed points $a$ (attractive) and $b$ (repulsive) in $\T$, 
and that the absolute value $|u|$ can be extended continuously to $\D \cup \{a,b\}$.
Then $uC_{\varphi}$ has the shadowing property if and only if one of the following conditions is true:
\begin{itemize}
\item [\rm (a)] $\displaystyle \frac{|u|(a)}{\varphi'(a)^{\gamma}} < \frac{|u|(b)}{\varphi'(b)^{\gamma}}$\,,
  and $|u|(a) > \varphi'(a)^{\gamma}$ or $ |u|(b) < \varphi'(b)^{\gamma}$.

\item [\rm (b)] $\displaystyle \frac{|u|(a)}{\varphi'(a)^{\gamma}} > \frac{|u|(b)}{\varphi'(b)^{\gamma}}$\,, 
  $|u|(a) \neq \varphi'(a)^{\gamma}$ and $|u|(b) \neq \varphi'(b) {\gamma}$.

\item [\rm (c)] $\displaystyle \frac{|u|(a)}{\varphi'(a)^{\gamma}} = \frac{|u|(b)}{\varphi'(b)^{\gamma}} \neq 1$ .
\end{itemize}
In the cases 
$\displaystyle \frac{|u|(a)}{\varphi'(a)^{\gamma}} < \frac{|u|(b)}{\varphi'(b)^{\gamma}}$
or
$\displaystyle \frac{|u|(a)}{\varphi'(a)^{\gamma}} = \frac{|u|(b)}{\varphi'(b)^{\gamma}}$, 
the operator $uC_{\varphi}$ has the shadowing property if and only if it is a hyperbolic operator.
\end{theorem}

\begin{proof}
Denote the weighted composition operator $uC_\varphi$ on $X$ simply by $T$.
Since $\partial \sigma(T) \subset \sigma_s(T)$, \cite[Theorem~8.1]{Oliva} implies that
\begin{itemize}
\item [(a)] If $\displaystyle \frac{|u|(a)}{\varphi'(a)^{\gamma}} < \frac{|u|(b)}{\varphi'(b)^{\gamma}}$, then 
  $\sigma_s(T) = \sigma(T) = \Big\{\lambda \in \C : \displaystyle \frac{|u|(a)}{\varphi'(a)^{\gamma}} \leq |\lambda| 
  \leq \frac{|u|(b)}{\varphi'(b)^{\gamma}}\Big\}$.

\item [(b)] If $\displaystyle \frac{|u|(a)}{\varphi'(a)^{\gamma}} > \frac{|u|(b)}{\varphi'(b)^{\gamma}}$, then 
  $\sigma_s(T) = \partial \sigma(T) = \Big\{\lambda \in \C : \displaystyle |\lambda| = \frac{|u|(a)}{\varphi'(a)^{\gamma}} 
  \mbox{ or } |\lambda| = \frac{|u|(b)}{\varphi'(b)^{\gamma}}\Big\}$ .

\item [(c)] If $\displaystyle \frac{|u|(a)}{\varphi'(a)^{\gamma}} = \frac{|u|(b)}{\varphi'(b)^{\gamma}}$, then 
  $\sigma_s(T) = \sigma(T) = \partial \sigma(T) = \Big\{\lambda \in \C : \displaystyle |\lambda| = \frac{|u|(a)}{\varphi'(a)^{\gamma}}\Big\}$.
\end{itemize}
Therefore, the result follows from the spectral characterization of the shadowing property.
\end{proof}

\begin{corollary} 
Let $X$ be a $\gamma$-space for some $\gamma > 0$ and let $C_{\varphi}$ be an invertible weighted composition operator on $X$. 
If $\varphi$ is a hyperbolic automorphism, then $C_{\varphi}$ has the shadowing property and is not a hyperbolic operator.
\end{corollary}

\begin {proof}
If $\varphi$ is a hyperbolic automorphism with fixed points $a$ (attractive) and $b$ (repulsive) in $\T$, 
then $0 < \varphi'(a) < 1$ and $\varphi'(b) = \frac{1}{\varphi'(a)}$\,; therefore, 
\[
\frac{1}{\varphi'(a)^{\gamma}} >\frac{1}{\varphi'(b)^{\gamma}}\,, \ \ 
1 \neq \varphi'(a)^{\gamma} \ \text{ and } \ 1 \neq \varphi'(b) ^{\gamma}.
\]
By case (b) of the theorem above, we conclude that $C_{\varphi}$ has the shadowing property.

On the other hand, by \cite[Theorem~8.1]{Oliva},
\[
\sigma(T) = \Big\{\lambda \in \C : \frac{1}{\varphi'(b)^{\gamma}} \leq |\lambda| \leq \frac{1}{\varphi'(a)^{\gamma}}\Big\} \supset \T.
\]
Thus, the operator $C_{\varphi}$ is not hyperbolic.
\end{proof}


\subsection{Shadowing for invertible weighted composition operators with parabolic or elliptic symbol}

Throughout this subsection, $X$ denotes a Banach space of holomorphic functions on the unit disk $\D$
which contains the constant functions and satisfies the following conditions:
\begin{itemize}
\item [(C1)] There exist positive constants $c$ and $s$ such that
\[
|f (z)| \leq c\, \|f\| (1-|z|^2)^{-s} \ \text{ for all } f \in X \text{ and } z \in \D,
\]
and, for each $z \in \D$, there exists $f_z \in X$ with $\|f_z\| \leq 1$ such that $f_z(z)(1 - |z|^2)^{s} =1$.
\item [(C2)] There exist positive constants $c$ and $s$ such that $\|C_{\varphi}\| \leq c\, (1 -|\varphi(0)|^2)^{-s}$
for all $\varphi \in Aut(\D)$.
\item [(C3)] For each $u \in H^{\infty}(\D)$, we have $\|M_u	\| \leq	\|u\|_{\infty}$.
\end{itemize}

The Hardy spaces $H^p(\D)$ and the weighted Bergman spaces $A^p_{\sigma}(\D)$ 
are the most important examples of spaces that satisfy (C1)--(C3).

\begin{theorem} 
Let $uC_{\varphi}$ be an invertible weighted composition operator on $X$. 
If the automorphism $\varphi$ is parabolic, $a$ is the unique fixed point of $\varphi$ in $\T$, and $u \in A(\D)$,
then the following assertions are equivalent:
\begin{itemize}
\item [\rm (i)] $uC_{\varphi}$ has the shadowing property;
\item [\rm (ii)] $uC_{\varphi}$ is hyperbolic;
\item [\rm (iii)] $|u|(a) \neq 1$.
\end{itemize}
\end{theorem}

\begin{proof}
Under the hypotheses of the theorem, \cite[Theorem~4.3]{HLNS} states that
\begin{equation}\label{Parabolic-Eq1}
\sigma(uC_\varphi) = \{\lambda \in \C : |\lambda| = |u|(a)\}.
\end{equation}
Since $\sigma_s(uC_\varphi)$ contains the boundary of $\sigma(uC_\varphi)$, (\ref{Parabolic-Eq1}) implies that
\begin{equation}\label{Parabolic-Eq2}
\sigma_s(uC_\varphi) = \sigma(uC_\varphi).
\end{equation}
The result follows from (\ref{Parabolic-Eq1}) and (\ref{Parabolic-Eq2}).
\end{proof}

\begin{corollary} 
Let $C_{\varphi}$ be an invertible composition operator on $X$. 
If the automorphism $\varphi$ is parabolic, then $C_{\varphi}$ does not have the shadowing property.
\end{corollary}

\begin{proof}
It follows from the previous theorem, since condition (iii) fails.
\end{proof}

\begin{theorem}
Let $uC_{\varphi}$ be an invertible weighted composition operator on $X$. 
If the automorphism $\varphi$ is elliptic with fixed point $a \in \D$, $\varphi^n \neq I_{\D}$ for all $n \in \N$, and $u \in A(\D)$,
then the following assertions are equivalent:
\begin{itemize}
\item [\rm (i)] $uC_{\varphi}$ has the shadowing property;
\item [\rm (ii)] $uC_{\varphi}$ is hyperbolic;
\item [\rm (iii)] $|u|(a) \neq 1$.
\end{itemize}
\end{theorem}

\begin{proof}
The proof is analogous to that of the previous theorem, 
but \cite[Theorem~4.14]{HLNS} should be used instead of \cite[Theorem~4.3]{HLNS}.
\end{proof}

\begin{corollary} 
Let $C_{\varphi}$ be an invertible composition operator on $X$. 
If the automorphism $\varphi$ is elliptic, then $C_{\varphi}$ does not have the shadowing property.
\end{corollary}

\begin{proof}
If $\varphi^n \neq I_{\D}$ for all $n \in \N$, then the result follows from the previous theorem, since condition (iii) fails.
If $\varphi^n = I_{\D}$ for a certain $n \in \N$, then $(C_\varphi)^n$ is the identity operator on $X$, 
which does not have the shadowing property. 
This implies that $C_\varphi$ also does not have the shadowing property.
\end{proof}


\section{Shadowing  property for weighted backward shifts}\label{Section7}


\subsection{Shadowing property for invertible bilateral weighted backward shifts} 

A complete characterization of the shadowing property for invertible bilateral weighted backward shifts on the classical Banach sequence spaces $\ell^p(\Z)$ ($1 \leq p < \infty$) and $c_0(\Z)$ was established in \cite[Theorem~18]{NBerAMes21}.
By combining Theorem~\ref{Spectral1} with the characterization of the approximate point spectrum of weighted shifts obtained by Ridge~\cite{R}, we now present a simple proof of this result for the case $1 < p < \infty$.

\begin{theorem}\label{ThmShadBw}
Let $X = \ell^p(\Z)$ $(1 < p < \infty)$ or $X = c_0(\Z)$.
Let $w = (w_n)_{n \in \Z}$ be a bounded sequence of scalars with $\inf_{n \in \Z} |w_n| > 0$ 
and consider the bilateral weighted backward shift
\[
B_w : (x_n)_{n \in \Z} \in X \mapsto (w_{n+1}x_{n+1})_{n \in \Z} \in X.
\] 
Then $B_w$ has the shadowing property if and only if one of the following conditions holds:
\begin{itemize}
\item [\rm (A)] $\displaystyle \lim_{n \to \infty} \sup_{k \in \Z} |w_k w_{k+1} \cdots w_{k+n}|^{\frac{1}{n}} < 1$.
\item [\rm (B)] $\displaystyle \lim_{n \to \infty} \inf_{k \in \Z} |w_k w_{k+1} \cdots w_{k+n}|^{\frac{1}{n}} > 1$.
\item [\rm (C)] $\displaystyle \lim_{n \to \infty} \sup_{k \in \N} |w_{-k} w_{-k-1} \cdots w_{-k-n}|^{\frac{1}{n}} < 1$ and 
  $\displaystyle \lim_{n \to \infty} \inf_{k \in \N} |w_k w_{k+1} \cdots w_{k+n}|^{\frac{1}{n}} > 1$.
\end{itemize}
\end{theorem}

\begin{proof}
Let $1 < q < \infty$ be the exponent conjugate to $p$, that is, $1/p + 1/q = 1$.
If $X = \ell^p(\Z)$ (resp.\ $X = c_0(\Z)$), 
then the adjoint operator $(B_w)^*$ of $B_w$ is the bilateral weighted forward shift $F_v$ on $\ell^q(\Z)$ (resp.\ $\ell^1(\Z)$) 
with weight sequence $v = (w_{n+1})_{n \in \Z}$.
Let 
\begin{align*}
i &= \lim_{n \to \infty} \inf_{k \in \Z} |w_k w_{k+1} \cdots w_{k+n-1}|^{\frac{1}{n}},\\
r &= \lim_{n \to \infty} \sup_{k \in \Z} |w_k w_{k+1} \cdots w_{k+n-1}|^{\frac{1}{n}},
\end{align*}
and
\begin{align*}
i^- &= \lim_{n \to \infty} \inf_{k \in \N} |w_{-k} w_{-k-1} \cdots w_{-k-n+1}|^{\frac{1}{n}},\\
r^- &= \lim_{n \to \infty} \sup_{k \in \N} |w_{-k} w_{-k-1} \cdots w_{-k-n+1}|^{\frac{1}{n}},\\
i^+ &= \lim_{n \to \infty} \inf_{k \in \N} |w_k w_{k+1} \cdots w_{k+n-1}|^{\frac{1}{n}},\\
r^+ &= \lim_{n \to \infty} \sup_{k \in \N} |w_k w_{k+1} \cdots w_{k+n-1}|^{\frac{1}{n}}.
\end{align*}
By arguing as in \cite[Theorem~3]{R} for the case $q = 2$ (see also \cite[Proposition~2.5]{BT}), it follows that:
\begin{itemize}
\item If $r^- < i^+$, then 
  $\sigma_a(F_v) = \{\lambda \in \C : i^- \leq |\lambda| \leq r^-\} \cup \{\lambda \in \C : i^+ \leq |\lambda| \leq r^+\}$.
\item If $r^- \geq i^+$, then $\sigma_a(F_v) = \{\lambda \in \C : i \leq |\lambda| \leq r\}$.
\end{itemize}
Since $\sigma_s(B_w) = \sigma_a(F_v)$, we see that the only possibilities for the intersection $\sigma_s(B_w) \cap \T$ to be empty are $r < 1$, $i > 1$, or $r^- < 1 < i^+$, which are equivalent to conditions (A), (B), and (C), respectively.
\end{proof}

The following dichotomy shows that the dynamics of any $B_w$ with the shadowing property is either simple or extremely chaotic.

\begin{corollary}\label{BWBS}
Let $X$, $w$ and $B_w$ be as in the previous theorem.
If $B_w$ has the shadowing property, then either $B_w$ is hyperbolic or $B_w$ satisfies the frequent hypercyclicity criterion.
\end{corollary}

The {\em frequent hypercyclicity criterion} \cite[Section~9.2]{GP} is a very strong property that simultaneously implies 
several types of chaotic behavior, including frequent hypercyclicity, topological mixing, Devaney chaos, dense distributional chaos, 
and dense mean Li-Yorke chaos (see \cite{BBMPe,BBPe,GP}).

\begin{proof}
Since $B_w$ has the shadowing property, Theorem~\ref{ThmShadBw} guarantees that 
one of the conditions (A), (B), or (C) must be satisfied.
Conditions (A) and (B) are exactly the cases in which $B_w$ is hyperbolic \cite[Remark~35(a)]{BCDMP}, 
and it was shown in \cite[Page~16]{NBerAMes-Arxiv} that condition (C) implies the frequent hypercyclicity criterion.
\end{proof}


\subsection{Positive shadowing property for unilateral weighted backward shifts} 

Let us now present a simple proof of \cite[Proposition~20]{NBerAMes-Arxiv} (see also \cite[Theorem~21]{BP}) for the case $1 < p < \infty$.

\begin{theorem}
Let $X = \ell^p(\N)$ $(1 < p < \infty)$ or $X = c_0(\N)$. 
Let $w = (w_n)_{n \in \N}$ be a bounded sequence of non-zero scalars and consider the unilateral weighted backward shift
\[
B_w : (x_1, x_2, \cdots) \in  X \mapsto (w_2 x_2, w_3 x_3, \cdots) \in  X.
\]
Then $B_w$ has the positive shadowing property if and only if one of the following conditions holds:
\begin{itemize}
\item [\rm (a)] $\displaystyle \lim_{n \to \infty} \sup_{k \in \N} |w_k w_{k+1} \cdots w_{k+n}|^{\frac{1}{n}} < 1$.
\item [\rm (b)] $\displaystyle \lim_{n \to \infty} \inf_{k \in \N } |w_k w_{k+1} \cdots w_{k+n}|^{\frac{1}{n}} > 1$.
\end{itemize}
\end{theorem}

\begin{proof}
Let $1 < q < \infty$ be the exponent conjugate to $p$.
If $X = \ell^p(\N)$ (resp.\ $X = c_0(\N)$), 
then the adjoint operator $(B_w)^*$ of $B_w$ is the unilateral weighted forward shift $F_v$ on $\ell^q(\N)$ (resp.\ $\ell^1(\N)$) 
with weight sequence $v = (w_{n+1})_{n \in \N}$.
Let 
\[
r = \lim_{n \to \infty} \sup_{k \in \N} |w_{k+1} w_{k+2} \cdots w_{k+n}|^{\frac{1}{n}} \ \ \text{ and } \ \ 
i = \lim_{n \to \infty} \inf_{k \in \N} |w_{k+1} w_{k+2} \cdots w_{k+n}|^{\frac{1}{n}}.
\]
By \cite[Proposrition~1.6.15]{LN},
\[
\sigma_a(F_v) = \{\lambda \in \C : i \leq |\lambda| \leq r\}.
\]
Since $\sigma_s(B_w) = \sigma_a(F_v)$, we see that the only possibilities for the intersection $\sigma_s(B_w) \cap \T$ to be empty 
are $r < 1$ or $i > 1$, which are equivalent to conditions (A) and (B), respectively.
\end{proof}

\begin{corollary}
Let $X$, $w$ and $B_w$ be as in the previous theorem.
If $B_w$ has the positive shadowing property, then either $B_w$ is hyperbolic or $B_w$ satisfies the frequent hypercyclicity criterion.
\end{corollary}

\begin{proof}
Condition (a) in the previous theorem is exactly the case in which $B_w$ is hyperbolic \cite[Remark~35(b)]{BCDMP}
and it was observed in \cite[Page~18]{NBerAMes-Arxiv} that condition (b) implies the frequent hypercyclicity criterion.
\end{proof}

\begin{theorem}
Let $X = \ell^p(\N)$ $(1 \leq p < \infty)$ or $X = c_0(\N)$. Consider the unilateral weighted backward shift
\[
B : (x_1, x_2, \cdots) \in X \mapsto (x_2, x_3, \cdots) \in  X.
\]
If $\phi : U \to \C$ is a non-constant holomorphic function on neighborhood of the closed unit disc, then the following assertions hold:
\begin{itemize}
\item [\rm (a)] $\phi(B)$ is hyperbolic if and only if $\phi(\ov{\D}) \cap \T = \emptyset$. 
\item [\rm (b)] $\phi(B)$ has the positive shadowing property if and only if $\phi(\T) \cap \T = \emptyset$.
\item [\rm (c)] $\phi(B)$ is Devaney chaotic if and only if $\phi(\D) \cap \T \neq \emptyset$.
\end{itemize}
\end{theorem}

\begin{proof}
(a): By the spectral mapping theorem, $\sigma(\phi(B)) = \phi(\sigma(B)) = \phi(\ov{\D})$, which gives (a).

\smallskip\noindent
(b):  By \cite[Theorem~2.48]{A}, $\sigma_s(\phi(B)) = \phi(\sigma_s(B)) = \phi(\T)$, which gives (b). 

\smallskip\noindent
(c): See \cite[Theorem~4.43]{GP}.
\end{proof}

\begin{remark}
As consequence of the previous theorem, we see that there exist operators of the form $\phi(B)$ on $\ell^p(\N)$ ($1 \leq p < \infty$)
and on $c_0(\N)$ that are Devaney chaotic and do not have the positive shadowing property.
However, every non-hyperbolic $\phi(B)$ with the positive shadowing property is Devaney chaotic.
In particular, the operator $I + \lambda B$ is Devaney chaotic and does not have the positive shadowing property 
if $0 < \lambda \leq 2$, and it is Devaney chaotic and has the positive shadowing property if $\lambda > 2$. 
\end{remark}

\begin{theorem} 
For an arbitrary unilateral weighted forward shift $T$ on $\ell^p(\N)$ $(1 \leq p < \infty)$, we have that:
\begin{itemize}
\item [\rm (a)] $T^*$ is hyperbolic if and only if $r(T) < 1$.
\item [\rm (b)] $T^*$ has the positive shadowing property and is not hyperbolic if and only if $i(T) > 1$.
\end{itemize}
\end{theorem}

\begin{proof}
(a): It follows from the fact that $\sigma (T) = \ov{\D}_{r(T)}$ \cite[Theorem~4.83]{A2}.

\smallskip\noindent
(b): Suppose that $i(T) > 1$. By (a), $T^*$ is not hyperbolic. Moreover, by \cite[Theorem~4.83]{A2},
\[
\sigma_s(T^*) = \sigma_a(T) = \{\lambda \in \C : i(T) \leq |\lambda| \leq r(T)\}.
\]
Hence, since $i(T) > 1$, $\sigma_s(T^*) \cap \T = \emptyset$ and $T^*$ has the positive shadowing property.

Conversely, suppose that $T^*$ has the positive shadowing property and is not hyperbolic.
By (a), $r(T) \geq 1$. Thus, if $i(T) \leq 1$, then $\sigma_s(T^*) \cap \T \neq \emptyset$ 
and $T^*$ does not have the positive shadowing property, a contradiction. Therefore, $i(T) > 1$.
\end{proof}


\section{Positive shadowing property for adjoints of multiplication operators}\label{Section8}

\begin{theorem} 
If $T \in L(X)$ is a non-invertible isometry and $\phi : U \to \C$ is a non-constant holomorphic function 
on a connected open neighborhood of the closed unit disc $\ov{\D}$, then the following assertions hold:
\begin{itemize}
\item [\rm (a)] $\phi(T)$ has the positive shadowing property if and only if $\phi(T)$ is hyperbolic 
  if and only if $\phi(\ov{\D}) \cap \T = \emptyset$.
\item [\rm (b)] $\phi(T)^*$ is hyperbolic if and only if $\phi(\ov{\D}) \cap \T = \emptyset$. 
\item [\rm (c)] $\phi(T)^*$ has the positive shadowing property if and only if $\phi(\T) \cap \T = \emptyset$.
\end{itemize}
\end{theorem}

\begin{proof}
(a): By \cite[Theorem~2.58(ii)]{A}, $\phi(T)$ has the SVEP.
Therefore, by Theorem~\ref{thm-SVEP}(a), $\phi(T)$ has the positive shadowing property if and only if it is hyperbolic.
Moreover, by \cite[Theorem~2.58(i)]{A}, $\sigma(\phi(T)) = \phi(\ov{\D})$. 
Thus, $\phi(T)$ is hyperbolic exactly when $\phi(\ov{\D}) \cap \T = \emptyset$.

\smallskip\noindent
(b): We have $\sigma(\phi(T)^*) = \sigma(\phi(T)) = \phi(\ov{\D})$, as we have already mentioned in (a).

\smallskip\noindent
(c): It follows from the fact that $\sigma_s(\phi (T)^*) = \sigma_a(\phi(T)) = \phi(\T)$, by \cite[Theorem~2.58(i)]{A}.
\end{proof}

Below, we present an application of the previous theorem to a more concrete situation.

\begin{corollary} 
Let $\phi : U \to \C$ be a non-constant holomorphic function on a connected open neighborhood of the closed unit disc $\ov{\D}$. 
Consider the {\em multiplication Toeplitz operator} $M_{\phi}$ on the Hardy space $H^2(\D)$, which is given by
\[
M_{\phi}(f) = \phi f \ \ \ \ \ (f \in H^2(\D)),
\]
the pointwise product of $\phi$ by $f$. The following assertions hold:
\begin{itemize}
\item [\rm (a)] $M_{\phi}$ has the positive shadowing property if and only if $M_{\phi}$ is hyperbolic 
  if and only if $\phi(\ov{\D}) \cap \T = \emptyset$.
\item [\rm (b)] $(M_{\phi})^*$ is hyperbolic if and only if $\phi(\ov{\D}) \cap \T = \emptyset$. 
\item [\rm (c)] $(M_{\phi})^*$ has the positive shadowing property if and only if $\phi(\T) \cap \T = \emptyset$.
\item [\rm (d)] $(M_{\phi})^*$ is Devaney chaotic if and only if $\phi(\D) \cap \T \neq \emptyset$.
\end{itemize}
\end{corollary}

\begin{proof}
Let $T$ be the operator on $H^2(\D)$ given by
\[
(Tf)(\lambda) = \lambda f(\lambda) \ \ \ \ \ (f \in H^2(\D), \lambda \in \D).
\]
The operator $T$ is a non-invertible isometry on $H^2(\D)$ and $\phi(T) = M_{\phi}$.
Therefore, items (a), (b) and (c) follow directly from the previous theorem.

Item (d) is a particular case of \cite[Theorem 4.42]{GP}.
\end{proof}

\begin{remark}
As consequence of the previous corollary, we see that there exist operators $(M_{\phi})^*$ on $H^2(\D)$
that are Devaney chaotic and do not have the positive shadowing property. 
However, every non-hyperbolic $(M_{\phi})^*$ with the positive shadowing property is Devaney chaotic.
\end{remark}

\begin{remark}
Note that similar results hold for multiplication Toeplitz operators with arbitrary bounded holomorphic symbols. 
In fact, for any $\phi \in H^{\infty}(\D)$, the surjective spectrum $\sigma_s((M_{\phi})^*)$ coincides with the essential range 
of the boundary function that is obtained by taking non-tangential limits of $\phi$ almost everywhere on the unit circle. 
\end{remark}


\section{Positive shadowing property for weighted composition operators on $L^2[0,1]$}\label{Section9}

According to \cite[Definition~2.1]{Po}, a pair of mappings $(\varphi,\psi)$ is called {\em $L^2$-admissible} 
if $\varphi : [0, 1] \to [0, 1]$ is an increasing differentiable function on $(0, 1)$ with only two fixed points $0$ and $1$,
and $\psi $  is a real positive-valued function on $[0, 1]$ such that
\[
\frac{\psi }{\sqrt{\varphi'}} \in L^{\infty}[0, 1].
\]
Moreover, if $\varphi (x) \geq x$ on $[0, 1]$, $(\varphi,\psi)$ is said to be {\em upper $L^2$-admissible},
and if $\varphi (x) \leq x$ on $[0, 1]$, $(\varphi,\psi)$ is said to be {\em lower $L^2$-admissible}.

It was proven in \cite[Proposition~2.4]{Po} that if the pair $(\varphi,\psi)$ is $L^2$-admissible,
then the weighted composition operator with symbol $\varphi$ and weight $\psi$ on $L^2[0, 1]$, 
\[
\psi C_{\varphi} : f \in L^2[0, 1] \mapsto \psi \cdot (f \circ \varphi) \in L^2[0, 1],
\]
is a well-defined bounded linear operator.

\begin{theorem} 
Let $(\varphi,\psi)$ be an upper $L^2$-admissible pair. 
Suppose that $g : [0, 1] \to \R$, defined by $g(t) = \frac{\psi}{\sqrt{\varphi'}}$, 
has a right-hand limit at $0$ denoted by $a$ and a left-hand limit at $1$ denoted by $b$. 
For $a \neq b$, the following assertions hold:
\begin{itemize}
\item [\rm (a)] If $a < b$, then $\psi C_{\varphi}$ has the positive shadowing property if and only if $a \neq 1$ and $b \neq 1$. 
\item [\rm (b)] If $a > b$, then $\psi C_{\varphi}$ has the positive shadowing property if and only if 
  $\psi C_{\varphi}$ is hyperbolic if and only if $a < 1$ or $b > 1$.  
\end{itemize}
\end{theorem}

\begin{proof}
(a): By \cite[Theorem~2.5]{Po}, $\sigma(\psi C_{\varphi}) = \{\lambda \in \C :  a \leq |\lambda| \leq b\}$
and, for $a < |\lambda| < b$, the adjoint $(\psi C_{\varphi} - \lambda I)^*$ is bounded below, 
which is equivalent to $\psi C_{\varphi} - \lambda I$ being surjective. Thus, 
\[
\sigma_s(\psi C_{\varphi}) = \{\lambda \in \C : |\lambda| = a \mbox{ or } |\lambda| = b\}.
\]
Therefore, $\psi C_{\varphi}$ has the positive shadowing property if and only if $a \neq 1$ and $b \neq 1$.

\smallskip\noindent
(b): In this case, it follows from \cite[Theorem~2.5]{Po} that
\[
\sigma(\psi C _{\varphi}) = \sigma_s (\psi C _{\varphi}) = \{\lambda \in \C : b \leq |\lambda| \leq a\}.
\]
Hence, $\psi C _{\varphi}$ has the positive shadowing property if and only if it is hyperbolic if and only if $a < 1$ or $b > 1$.  
\end{proof}

\begin{theorem} 
Let $(\varphi,\psi)$ be a lower $L^2$-admissible pair. 
Suppose that $g : [0, 1] \to \R$, defined by $g(t) = \frac{\psi}{\sqrt{\varphi'}}$, 
has a right-hand limit at $0$ denoted by $a$ and a left-hand limit at $1$ denoted by $b$. 
For $a \neq b$, the following assertions hold:
\begin{itemize}
\item [\rm (a)] If $a < b$, then $\psi C_{\varphi}$ has the positive shadowing property if and only if 
  $(\psi C_{\varphi})$ is hyperbolic if and only if $a > 1$ or $b < 1$. 
\item [\rm (b)] If $a > b$, then $\psi C_{\varphi}$ has the positive shadowing property if and only if $a \neq 1$ and $b \neq 1$.  
\end{itemize}
\end{theorem}

\begin {proof}
(a): It follows from \cite[Theorem 2.7]{Po} that
\[
\sigma(\psi C_{\varphi}) = \sigma_s(\psi C_{\varphi}) = \{\lambda \in \C :  a \leq |\lambda| \leq b\}.
\]
Thus, $(\psi C_{\varphi})$ has the positive shadowing property if and only if it is hyperbolic if and only if $a > 1$ or $b < 1$. 

\smallskip\noindent
(b): In this case, \cite[Theorem 2.7]{Po} tell us that $\sigma(\psi C_{\varphi}) = \{\lambda \in \C :  b \leq |\lambda| \leq a\}$
and, for $b < |\lambda| < a$, $(\psi C_{\varphi} - \lambda I)^*$ is bounded below, that is, $\psi C_{\varphi} - \lambda I$ is surjective.
Thus, 
\[
\sigma_s(\psi C_{\varphi}) = \{\lambda \in \C : |\lambda| = a  \mbox{ or } |\lambda | = b\},
\]
which implies the desired result.  
\end{proof}


\section{Shadowing property for weighted composition operators on $L^p(\mu)$}\label{Section10}

Throughout this section, we fix $p \in [1,\infty)$.
Recall \cite[Definition 13]{BBP2} that a {\em weighted $p$-measurable system} is a $5$-tuple $(X,\fM,\mu,f,w)$ such that 
the following conditions hold:
$(X,\fM,\mu)$ is a $\sigma$-finite measure space with $\mu(X) > 0$;
$f : X \to X$ is a non-singular bimeasurable bijective map, 
where $f$ {\em non-singular} means that $\mu(f^{-1}(B)) = 0$ if and only if $\mu(B) = 0$;
and $w : X \to \K$ is a bounded measurable map for which there is a constant $c > 0$ such that
\[
\int_B |w|^p d\mu \leq c\, \mu(f(B)) \ \text{ for all } B \in \fM.
\]
This implies that the weighted composition operator 
\[
C_{w,f}(\varphi) = w \cdot (\varphi \circ f)
\]
is a well-defined bounded linear operator on $L^p(\mu)$.
Note that, in this section, we follow the notation of \cite{BBP,BBP2} 
and use $C_{w,f}$ instead of $wC_f$ for the weighted composition operator described above.
If $w$ is the constant function $1$, then we obtain the concept of a {\em measurable system} $(X,\fM,\mu,f)$.
Recall that such a system (or the map $f$) is said to be {\em dissipative} if there is a measurable set $W$ (called a {\em wandering set}) 
such that the sets $f^n(W)$, $n \in \Z$, are pairwise disjoint and $X = \bigcup_{n \in \Z} f^n(W)$. 

Recall \cite[Definition~14]{BBP2} that a weighted $p$-measurable system $(X,\fM,\mu,f,w)$ is said to be {\em invertible} if 
$\wt{w} = \frac{1}{w \circ f^{-1}}$ is a well-defined bounded measurable map on $X$ and there is a constant $\tilde{c} > 0$ satisfying
\[
\int_B |\wt{w}|^p d\mu \leq \tilde{c}\, \mu(f^{-1}(B)) \ \text{ for all } B \in \fM.
\]
In this case, $C_{w,f}$ is an invertible operator on $L^p(\mu)$ and $(C_{w,f})^{-1} = C_{\wt{w},f^{-1}}$.
Moreover, we consider the following bilateral sequence of positive measures on $(X,\fM)$:
\[
\mu_n(B) = \int_B |w^{(n)}|^p d\mu \ \ \ \ \ (B \in \fM, n \in \Z),
\]
where $w^{(0)} = 1$, $w^{(1)} = w$, $w^{(-1)} = \wt{w}$,
\[
w^{(n)} = (w \circ f^{n-1}) \cdots (w \circ f) \cdot w \ \text{ and } \ 
w^{(-n)} = (\wt{w} \circ f^{-(n-1)}) \cdots (\wt{w} \circ f^{-1}) \cdot \wt{w} \ \ \ (n \geq 2).
\]

The dynamics of (weighted) composition operators on $L^p(\mu)$ spaces has been extensively investigated 
by various authors over the last few years. 
We refer the reader to \cite{BBP,BBP2,BV,BDDP,DDM2}, where additional references can be found. 
The following result presents sufficient conditions for the shadowing property in the present context.

\begin{theorem}\label{ShadLp}
Let $(X,\fM,\mu,f,w)$ be a weighted $p$-measurable system.
\begin{itemize}
\item [\rm (a)] If 
  $\displaystyle \lim_{n \to \infty} \Big(\sup_{0 < \mu(B) < \infty} \frac{\mu_n(f^{-n}(B))}{\mu(B)}\Big)^{\frac{1}{n}} < 1$,
  then $C_{w,f}$ has the positive shadowing property.
\item [\rm (b)] If $(X,\fM,\mu,f,w)$ is invertible and one of the following conditions holds:
  \begin{itemize}
  \item [\rm (b1)] $\displaystyle \lim_{n \to \infty} \Big(\sup_{0 < \mu(B) < \infty} \frac{\mu_n(f^{-n}(B))}{\mu(B)}\Big)^{\frac{1}{n}} < 1$, or
  \item [\rm (b2)] $\displaystyle \lim_{n \to \infty} \Big(\sup_{0 < \mu(B) < \infty} \frac{\mu_{-n}(f^{n}(B))}{\mu(B)}\Big)^{\frac{1}{n}} < 1$, or
  \item [\rm (b3)] $(X,\fM,\mu,f,w)$ is dissipative with wandering set $W$, 
    \[
    \lim_{n \to \infty} \Big(\sup_{B\, \subset\, \bigcup_{k=1}^{\infty}f^{-k}(W), \; 0 < \mu(B) < \infty} 
      \frac{\mu_n(f^{-n}(B))}{\mu(B)}\Big)^{\frac{1}{n}} < 1
    \]
    and
    \[ 
    \lim_{n \to \infty} \Big(\sup_{B\, \subset\, \bigcup_{k=0} ^{\infty}f^{k}(W), \; 0 < \mu(B) < \infty} 
      \frac{\mu_{-n}(f^{n}(B))}{\mu(B)}\Big)^{\frac{1}{n}} < 1,
    \]
  \end{itemize}
then $C_{w,f}$ has the shadowing property.
\end{itemize}
\end{theorem}

\begin{proof}
(a): By \cite[Theorem~6]{BBP},
\[
\|(C_{w,f})^n\| = \sup_{0 < \mu(B) < \infty} \Big(\frac{\mu_n(f^{-n}(B))}{\mu(B)}\Big)^\frac{1}{p} \ \ \text{ for all } n \in \N.
\]
Thus,
\[
r(C_{w,f}) = \lim_{n \to \infty} \Big(\sup_{0 < \mu(B) < \infty} \Big(\frac{\mu_n(f^{-n}(B))}{\mu(B)}\Big)^\frac{1}{p}\Big)^{\frac{1}{n}} < 1,
\]
by the hypothesis.
This shows that $C_{w,f}$ is hyperbolic and, therefore, has the positive shadowing property (Corollary~\ref{GH-PSP}).

\smallskip\noindent
(b): We have already seen in item (a) that condition (b1) implies that $C_{w,f}$ is hyperbolic with $r(C_{w,f}) < 1$.
Similarly, by \cite[Theorem~6]{BBP}, condition (b2) implies that $C_{w,f}$ is hyperbolic with
\[
r((C_{w,f})^{-1}) = \lim_{n \to \infty} \Big(\sup_{0 < \mu(B) < \infty} \Big(\frac{\mu_{-n}(f^{n}(B))}{\mu(B)}\Big)^\frac{1}{p}\Big)^{\frac{1}{n}}< 1.
\]
Suppose that condition (b3) holds.
For each $\varphi \in L^p(\mu)$, define
\begin{align*}
\varphi_+(x) &= \left\{
  \begin{array}{cl}
    0 & \text{if } x \in \bigcup_{n=0}^{\infty}f^n(W)\\
    \varphi(x) & \text{otherwise,}
  \end{array}
	     \right.
\end{align*}
\begin{align*}
\varphi_-(x) &= \left\{
  \begin{array}{cl}
   0 & \text{if } x \in \bigcup_{n=1}^{\infty} f^{-n}(W)\\
    \varphi(x) & \text{otherwise.}
  \end{array}
	     \right.
\end{align*}
Let $L_+ = \{\varphi_+ : \varphi \in L^p(\mu)\}$ and $L_- = \{\varphi_- : \varphi \in L^p(\mu)\}$. 
Clearly,
\[
L^p(\mu) =L_+\oplus L_-, \ \ \ C_{w,f}(L_+) \subset L_+ \ \ \text{ and } \ \ (C_{w,f})^{-1}(L_-) \subset L_-.
\]
Furthermore, for each $n \in \N$,
\[
\|((C_{w,f})|_{L_+})^n\| = \sup_{B\, \subset\, \bigcup_{k=1}^{\infty} f^{-k}(W), \; 0 < \mu(B) < \infty} 
  \Big(\frac{\mu_n(f^{-n}(B))}{\mu(B)}\Big)^\frac{1}{p},
\]
\[
\|(((C_{w,f})^{-1})|_{L_-})^n\| = \sup_{B\, \subset\, \bigcup_{k=0}^{\infty} f^{k}(W), \; 0 < \mu(B) < \infty} 
  \Big(\frac{\mu_{-n}(f^{n}(B))}{\mu(B)}\Big)^\frac{1}{p}.
\]
Thus, our assumptions in (b3) imply that $r((C_{w,f})|_{L_+}) < 1$ and $r(((C_{w,f})^{-1})|_{L_-}) < 1$,
showing that $C_{w,f}$ is generalized hyperbolic.
In any case, we conclude that $C_{w,f}$ has the shadowing property (Corollary~\ref{GH-SP}).
\end{proof}

\begin{remark}
If $X = \Z$, $\fM = \cP(\Z)$ (the power set of $\Z$), $\mu$ is the counting measure on $\fM$,
$f$ is the $+1$ map, and $w = (w_n)_{n \in \Z}$ is a bounded sequence of scalars with $\inf_{n \in \Z} |w_n| > 0$,
then $L^p(\mu) = l^p(\Z)$ and $C_{w,f}$ is the invertible bilateral weighted backward shift
\[
B_w : (x_n)_{n \in \Z} \in \ell^p(\Z) \mapsto (w_{n}x_{n+1})_{n \in \Z} \in \ell^p(\Z).
\]
In this case, conditions (b1), (b2) and (b3) in Theorem~\ref{ShadLp} are equivalent to
conditions (A), (B) and (C) in Theorem~\ref{ThmShadBw}, which characterize the shadowing property for $B_w$.
This shows that, in this particular case, the sufficient conditions in Theorem~\ref{ShadLp}(b) are also necessary.
\end{remark}


\section*{Acknowledgement}

The first author was partially supported by CNPq – Project {\#}318556/2025-1, and by CAPES
– Finance Code 001. The first and the second authors were partially supported by Project
PID2022-139449NB-I00, funded by MCIN/AEI/10.13039/501100011033/FEDER, UE.


\end{document}